\documentclass[11pt]{article}
\usepackage{amssymb,amsfonts,amsmath,amstext,amsthm,amscd}

\newtheorem{thm}{Theorem}[section]
\newtheorem{prop}[thm]{Proposition}
\newtheorem{lemma}[thm]{Lemma}
\newtheorem{prob}[thm]{Problem}
\newtheorem{cor}[thm]{Corollary}

 \makeatletter
\newfont{\footsc}{cmcsc10 at 8truept}
\newfont{\footbf}{cmbx10 at 8truept}
\newfont{\footrm}{cmr10 at 10truept}
\title{The Strong Primitive Normal Basis Theorem}

\author {Stephen D. Cohen and Sophie Huczynska}

\date{}

\begin{document}
\maketitle

\newcommand{\Es} {E^{*}}
\newcommand{\Fs}{F^{*}}
\newcommand{\alp}{\alpha \in E}
\newcommand{\n}{n^*}
\newcommand{\w}{\omega}
\newcommand{\Q}{Q^*}
\newcommand{\ic}{\textit}
\numberwithin{equation}{section}
\def\mod{\mathop{\mathrm{mod}}}
\def\deg{\mathop{\mathrm{deg}}}

\renewcommand{\baselinestretch}{1.5}

\abstract {An element $\alpha$ of the extension $E$ of degree $n$
over the finite field $F=GF(q)$ is called \emph{free over $F$} if
$\{\alpha, \alpha^q, \ldots, \alpha^{q^{n- 1}}\}$ is a (normal)
basis of $E/F$.  The \emph{primitive normal basis theorem}, first
established in full by Lenstra and Schoof (1987), asserts that for
any such extension $E/F$, there exists an element $\alpha \in E$
such that $\alpha$ is simultaneously primitive (i.e., generates the
multiplicative group of $E$) and free over $F$.  In this paper we
prove the following strengthening of this theorem: aside from five
specific extensions $E/F$, there exists an element $\alpha \in E$
such that both $\alpha$ and $\alpha^{-1}$ are simultaneously
primitive and free over $F$.} 

\section{Introduction}\label{sec:1}
Given $q$, a power of a prime $p$, denote by $F$ the finite field
$GF(q)$ of order $q$, and by $E$ its extension $GF(q^n)$ of degree
$n$.  A \ic{primitive element} of $E$ is a generator of the cyclic
group $E^*$.  Additively, too, the extension $E$ is cyclic when
viewed as an $FG$-module, $G$ being the Galois group of $E$ over
$F$.  The classical form of this result - \ic{the normal basis
theorem} - is stated as follows:

\begin{thm}[Normal Basis Theorem]
There exists an element $\alpha \in E$ (an additive generator)
whose conjugates $\{ \alpha, \alpha^q, \dots, \alpha^{q^{n-1}} \}$
form a basis of $E$ over $F$.
\end{thm}

Such an element $\alpha$ is a \ic{free (or normal) element} of $E$
over $F$, and a basis of this kind is a \ic{normal basis} over
$F$. The key existence result linking additive and multiplicative
structure is \ic{the primitive normal basis theorem}:

\begin{thm}[Primitive Normal Basis Theorem]\label{PNBT}
For every prime power $q$ and $n \in \mathbb{N}$, there exists
$\alpha \in E$, simultaneously primitive and free over $F$.
Equivalently, there exists a \ic{primitive normal basis} over $F$,
all of whose members are primitive and free.
\end{thm}

Existence of such a basis for every extension was first proved by
Lenstra and Schoof \cite{LeSc87}, completing work by Carlitz
\cite{Car1}, \cite{Car2}, and Davenport \cite{Dav68}. A
computer-free proof of this result was produced by Cohen and
Huczynska \cite{CoHu03}.  The key to the transition to the more
theoretical and less computational approach realised in
\cite{CoHu03} was the introduction of sieving techniques (cf.
Section \ref{Sieve}, below). The question arises as to whether a yet
stronger existence theorem concerning primitive and free elements
can be proved unconditionally (or with very few exceptions) by means
of such techniques.  In this paper, we consider the following
natural problem, first suggested to us by  Robin J. Chapman (Exeter)
(to whom we are grateful).

\begin{prob}[PFF problem] \label{PFF} Given a finite extension $E/F$
of Galois fields, does there exist a primitive element $\alpha$ of
$E$, free over $F$, such that its reciprocal $\alpha^{-1} \in E$ is
also primitive and free over $F$? If so, then the pair $(q,n)$
corresponding to $E/F$ is called a PFF-pair.
\end{prob}

Observe that, for $\alpha \in E$, $\alpha$ is a primitive element
of $E$ if and only if $\alpha^{-1}$ is primitive; hence the four
conditions in Problem \ref{PFF} effectively reduce to three
($\alpha$ \emph{p}rimitive and \emph{f}ree, $\alpha^{-1}$
\emph{f}ree).

In this paper, we solve this problem completely: the answer is in
the affirmative except for a small number of listed exceptions. We
obtain the following strengthening of the Primitive Normal Basis
Theorem.

\begin{thm}[Strong Primitive Normal Basis Theorem]\label{SPNBT}
For every prime power $q$ and $n \in \mathbb{N}$, there exists a
primitive element $\alpha$ of $E$, free over $F$, such that its
reciprocal $\alpha^{-1} \in E$ is also primitive and free over
$F$, unless the pair $(q,n)$ is one of
\[(2,3),\;(2,4),\;(3,4),\; (4,3),\;(5,4).\]
\end{thm}

 Towards Theorem $\ref{SPNBT}$, Tian and Qi $\cite{TQ06}$
have given a proof provided $n \geq 32$ (when there are no
exceptions). They use an elaboration of the method of Lenstra and
Schoof \cite{LeSc87} but do not employ any of the sieving techniques
that are a feature of the present article and appear to be necessary
for completion, particularly for small values of $n$. Moreover,
because of the demanding nature of the PFF condition, fields of
smallest cardinality require individual treatment.  Our
consideration of the general problem therefore takes place in the
setting where $q \geq 5$ (even here, special care is needed for
$q=5$ and $7$), and we deal with the case $2 \leq q \leq 4$ in
Section \ref{vsmall} ``Very small fields". In what follows, all
non-trivial computation is performed using MAPLE (Version 10). Aside
from the five genuine exceptions listed in Theorem $\ref{SPNBT}$,
there are $35$ pairs $(q,n)$ (with $q\leq 13, n \leq 16$) for which
verification is by direct construction of a PFF polynomial:
otherwise, the proof is purely theoretical.

\section{Reductions}\label{sec:2}
In this section, we formulate the basic theory and perform some
reductions to the problem.  As much as possible, we aim to make
this account self-contained.

We begin by extending the notions of primitivity and free-ness. Let
$w\in\Es$. Then $w$ is a primitive element of $E$ if and only if $w$
has multiplicative order $q^n-1$, i.e., $w=v^d$ ($w \in E$) implies
$(d,q^n-1)=1$.  We extend this concept as follows: for any divisor
$m$ of $q^n-1$, we say that $w \in E^*$ is $m$-\ic{free}, if $w=v^d$
(where $v \in E$ and $d|m$) implies $d=1$. Thus $w \in E^*$ is
$m$-free if and only if $w$ is not an $l$th power for all primes $l$
dividing $m$. It follows that $w$ is $m$-free if and only if it is
$m_0$-free, where $m_0$ is the radical of $m$, i.e., the product of
its distinct prime factors.   In the context of the PFF problem,
observe that $w$ is $m$-free if and only if $w^{-1}$ is $m$-free
since, if $w^{-1}=v^k$ for some $k|m$ and $v \in \Es$, then
$w=(v^{-1})^k$ and $v^{-1} \in \Es$.

For $w\in E$, the $F$-order of $w$ is defined to be the monic
divisor $g$ (over $F$) of $x^n-1$ of minimal degree such that
$g^{\sigma}(w)=0$ ($g^{\sigma}$ is the polynomial obtained from $g$
by replacing each $x^i$ with $x^{q^i}$). Clearly, $w$ is free if and
only if the $F$-order of $w$ is $x^n-1$.  If $w \in E$ has $F$-order
$g$, then $w=h^{\sigma}(v)$ for some $v \in E$, where
$h=\frac{x^n-1}{g}$. Let $M$ be an $F$-divisor of $x^n-1$.  If
$w=h^{\sigma}(v)$ (where $v \in E$, $h$ is an $F$-divisor of $M$)
implies $h=1$ we say that $w$ is  $M$-\ic{free in} $E$. Again, $M$
may be replaced by its radical.  An important instance of this
occurs when $n$ is divisible by the characteristic $p$, say
$n=p^b\n$ (where $p \nmid \n$), in which event $w$ is $x^n-1$-free
if and only if it is $x^{\n} -1$-free. (The expansion of $n=p^b\n$,
as above, will be assumed throughout.)

We remark that, in the sequel, most arguments concerning divisors of
a given integer divisor of $q^n-1$ or polynomial divisors of a given
factor of $x^n-1$ depend only on the appropriate radicals  so that
the divisors may be assumed to be square-free.  To avoid awkward
qualifications to  these arguments, the reader is requested
throughout to interpret all relevant statements accordingly.

We make the following observation.

\begin{lemma}
Let $x^d-w$ be an $F$-divisor of $x^n-1$ ($w \in \Fs$, $d|n$).
Then, for $\alpha \in \Es$,
\[ (x^d-w)^{\sigma}(\alpha)=0 \Leftrightarrow (x^d-
w^{-1})^{\sigma}({\alpha}^{-1})=0. \] In particular, if $w \in
\Es$ has $F$-order $x+1$ or $x-1$, then so does $w^{-1}$.
\end{lemma}

If $n=2$ and $w\in\Es$ is primitive, then neither $w$ nor $w^{-1}$
can have $F$-order $ x\pm 1$ and so both are free over $F$.
Henceforth, we assume $n \geq 3$.

\begin{lemma} \label{n prime}
Let $n$ ($\geq 5$) be prime.  Suppose that $q$ is such that $p \nmid
n$ and $q \,(\mod \, n)$ is a multiplicative generator of the cyclic
group $(\mathbb{Z}/n \mathbb{Z})^*$.  Then $(q,n)$ is a PFF pair.
\end{lemma}
\begin{proof}  Under the given circumstances, $(n,q^i-1)=1$ for
$i=1, \ldots, n-2$ and $n$ for $i=n-1$; so $x^n-1$ factorizes into
irreducibles over $F$ as $(x-1)(x^{n-1}+x^{n-2}+ \ldots +x+1)$. By
Theorem 1.1 of \cite{Coh200}, there exists a primitive element $w
\in E$ such that its trace over $F$, $T(w) \neq 0$ and, similarly,
$T(w^{-1}) \neq 0$, i.e. neither $w$ nor $w^{-1}$ has $F$-order
$x^{n-1}+x^{n-2}+ \ldots +x+1$.  Since $w$ is primitive, neither $w$
nor $w^{-1}$ can have $F$-order $x-1$.
\end{proof}

Observe that Lemma \ref{n prime} applies to $\phi(n-1)$ of the $n$
possible congruence classes for values of $q$.  The next result
demonstrates the application of the lemma to some small values of
$n$.

\begin{lemma}\label{EarlyEx}
For the following values of $q$ and $n$, the pair $(q,n)$ is a PFF
pair:
\end{lemma}
\begin{itemize}
\item[(i)] $n=5$; $q \equiv 2$ or $3 \, (\mod\, 5)$.

\item[(ii)] $n=7$; $q \equiv 3$ or $5 \, (\mod\, 7)$.

\item[(iii)] $n=11$; $q \equiv 2$, $6$, $7$ or $8 \, (\mod\, 11)$.
\end{itemize}

For any $m|\,q^n-1$, and $g, \, h|\,x^n-1$, denote by $N(m,g,h)$
the number of non-zero elements $w \in E$ such that $w$ is
$m$-free and $g$-free, and $w^{-1}$ is $h$-free (note that
$w^{-1}$ is automatically $m$-free). As a consequence of the
earlier discussion, we may replace $m$, $g$ or $h$ by their
radicals at any time. To solve the PFF problem it would suffice to
show that $N(q^n-1,x^n-1,x^n-1)$ is positive, for every pair
$(q,n)$; however, it is useful to refine this requirement.

For a given pair $(q,n)$, define   $Q:= Q(q,n)$ to be (the radical
of) $\frac{q^n-1}{(q-1)\gcd(n,q-1)}$.  As in  \cite{LeSc87} and
\cite{Coh200}, we now demonstrate that  $q^n-1$ may be replaced by
$Q$, i.e. it suffices to show that $N(Q,x^n-1,x^n-1)$ is positive.
The following lemma, analogous to Lemma 2.1 of \cite{CoHu03}, makes
this relationship explicit.

\begin{lemma}\label{Q}
For any pair (q,n), $$N(Q,x^n-1,x^n-1)=
\frac{R}{\phi(R)}N(q^n-1,x^n-1,x^n-1),$$ where $\phi$ denotes
Euler's function, and $R$ is the greatest divisor of $q^n-1$
co-prime to $Q$.
\end{lemma}
\begin{proof} Let $\Q:=(q^n-1)/R$: then $\Q$ is the greatest divisor of
$q^n-1$ whose prime factors are those of $Q$. Moreover, $Q |\,\Q,
\; R |\,(q-1)(n,q-1) | \, (q-1)^2$, and $(R,\Q)=1$. In particular,
if $\gamma (\in \Es)$ is an $R$-th root of unity, then
$c:=\gamma^{q-1} \in \Fs$, and $\gamma^{q^i}=c^i\gamma$, for every
$i$. It follows  that, if $\alp$  and $\gamma$ is any $R$-th root
of unity, then $\alpha$ is $x^n-1$ free if and only if
$\gamma\alpha$ is  $x^n-1$ free. (Indeed, for any $k$, with $0
\leq k < n$,
 \[\sum_{i=0}^k a_i(\gamma\alpha)^{q^i}=0 \iff \sum_{i=0}^k
 a_ic^i\alpha^{q^i}=0, \quad a_0, \ldots, a_k, \,c \in F. ) \]

Now, any element $\alpha \in E^*$ can be expressed uniquely as the
product of a $\Q$-th root of unity $\alpha_0$ and an $R$-th root of
unity (in $\Es$).  By the above, if $\alpha$ is $Q$-free and  both
$\alpha$ and $\alpha^{-1}$ are $x^n-1$-free, then $\gamma\alpha_0$
is also $Q$-free with $\gamma\alpha_0$ and its inverse both
$x^n-1$-free, for any $R$-th root of unity $\gamma$.  If in fact
$\alpha$ is primitive, then $\alpha = \gamma\alpha_0$, for some
primitive $R$-th root of unity, $\gamma$.
\end{proof}

The following result will prove useful.
\begin{lemma} \label{Red2}
\begin{itemize}
\item[(i)] Assume $n=4$ and $q \equiv 3 \, (\mod\, 4)$. Then
$N(Q,x^4-1,x^4-1)=N(Q,x^2-1,x^2-1)$.

\item[(ii)] Assume $n=3$ and $q \equiv 2 \,(\mod\, 3)$. Then
$N(Q,x^3-1,x^3-1)=N(Q,x-1,x-1)$.
\end{itemize}
\end{lemma}
\begin{proof}  Take the case with  $n=4$, so that $x^2+1$ is
irreducible over $F$.  Suppose that $\alpha$ is $Q$-free and
$x^2-1$-free and $\alpha^{-1}$ is $x^2-1$-free, but $\alpha$ is not
$x^4-1$-free. Then $\alpha= \beta^{q^2} + \beta$, and hence
$\alpha^{q^2}= \alpha$, i.e., $\alpha^{q^2-1}=1$. Thus $\alpha$ is
contained in the quadratic extension of $F$ and so cannot be
primitive.  The same argument ensures that $\alpha^{-1}$ is also
$x^4-1$-free. The ``$n=3$" case is exactly analogous.
\end{proof}

\section{An expression for $N(m,g,h)$}\label{sec:3}

In this section, we employ character sums to obtain expressions,
and thence estimates, for the number of elements of the desired
type. We suppose throughout that $m|Q$ and $g, \, h|x^n-1$, where,
if desired, these can be assumed to be square-free.  We begin by
establishing characteristic functions for those subsets of $E$
comprising elements that are $m$-free, $g$-free or $h$-free.

\bigskip

\noindent \textbf{I. } \textit {The set of  $w \in \Es$ that are
$m$-free.}

Let $\hat{\Es} (\cong \Es)$ denote the group of multiplicative
characters of $\Es$.  For any $d|Q$, we write $\eta_d$ for a
typical character in $\hat{\Es}$ of order $d$.  Thus $\eta_1$ is
the trivial character. Notice that, since $d|\frac{q^n-1}{q-1}$,
the restriction of $\eta_d$ to $\Fs$ is the trivial character
$\nu_1$ of $\hat{\Fs}$.

We employ the following notation for weighted sums (cf.
\cite{CoHu03}).  For $m|Q$, set
\[ \int\limits_{d|m} \eta_{d} :=
\sum_{d|m} \frac{\mu(d)}{\phi(d)} \sum_{(d)} \eta_{d},
\]
where $\phi$ and $\mu$ denote the functions of Euler and M\"{o}bius
respectively and the inner sum runs through all $\phi(d)$ characters
of order $d$. (Note that only square-free divisors $d$ have any
influence.) Then, according to a formula developed from one of
Vinogradov, the characteristic function for the subset of $m$-free
elements of $\Es$ is

\begin{equation} \label{Vin} \theta(m) \int_{d|m} \eta_{d}
(w),\quad w \in \Es, \end{equation}
 where $\theta(m):=
\frac{\phi(m)}{m} = \prod\limits_{l|m, \;l\, \mathrm{prime}}
(1-l^{-1})$. (In Vinogradov's original formula characterising
primitive roots of a prime $p$, (\ref{Vin}) holds with $m=p-1$.)

\bigskip
\noindent \textbf{II. }\textit{The set of  $w \in E$ that are
$g$-free or $h$-free over $F$.}

Let $\lambda$ be the canonical additive character of $F$. Thus,
for $x \in F$,
\[ \lambda(x)=\mathrm{exp}(2\pi i
T_{F,\mathbb{F}_p}(x)/p),
\]
where $p$ is the characteristic of $F$ and $T_{F,\mathbb{F}_p}$
denotes the trace function from $F$ to $\mathbb{F}_p$.

Now let $\chi$ be the canonical additive character on $E$; it is
simply the lift of $\lambda$ to $E$, ie. $\chi(w)=\lambda(T(w))$, $w
\in E$. For any (monic) $F$-divisor $D$ of $x^n-1$, a typical
character $\chi_D$ of $F$- order $D$ is one such that $\chi_D \circ
D^{\sigma}$ is the trivial character in $E$, and $D$ is minimal (in
terms of degree) with this property.  For any $\delta\in E$, let
$\chi_\delta$ be the character defined by
$\chi_{\delta}(w)=\chi(\delta w),w \in E$. Define the subset
$\Delta_D$ of $E$ as the set of $\delta$ for which $\chi_{\delta}$
has $F$-order $D$. So we may also write $\chi_{\delta_D}$ for
$\chi_D$, where $\delta_D \in \Delta_D$; moreover $\{
\chi_{\delta_D}:\delta_D \in \Delta_D \}$ is the set of all
characters of order $D$.  Note that $\Delta_D$ is invariant under
multiplication by $\Fs$, and that, if $D=1$, then $\delta_1 =0$ and
$\chi_D = \chi_0$, the trivial character.  There are $\Phi(D)$
characters $\chi_D$, where $\Phi$ is the Euler function on $F[x]$
($\Phi$ is multiplicative and is given by the formula $\Phi(D)=|D|
\prod\limits_{P|D} (1-|P|^{-1})$, where the product is over all
monic irreducible $F$-divisors of $D$ and
$|D|=q^{\mathrm{deg}(D)}$).

In analogy to \textbf{I}, for $g|x^n-1$, define
\[ \int\limits_{D|g} \chi_{\delta_D} := \sum\limits_{D|g}
\frac{\mu(D)}{\Phi(D)} \sum\limits_{\delta_D} \chi_{\delta_D}, \]
where $\mu$ is the M\"{o}bius function on $F[x]$ and the inner sum
runs through all $\Phi(D)$ elements $\delta_D$ of $\Delta_D$ (only
square-free $D$ matter). With the notation
$\Theta(g)=\frac{\Phi(g)}{|g|}$, the characteristic function of the
set of $g$-free elements of $E$ correspondingly takes the form
\[ \Theta(g) \int\limits_{D|g} \chi_{\delta_D}(w), \quad w \in E.
\]

Using these characteristic functions, we derive an expression for
$N(m,g,h)$ in terms of Kloosterman and Gauss sums on $E$ and $F$.
For any $\alpha$, $\beta \in E$ and any multiplicative character
$\eta \in \hat{\Es}$, we define the generalized Kloosterman sum
$K(\alpha,\beta;\eta)$ ($=K_{q,n}(\alpha,\beta;\eta)$) by
\begin{equation*}
K(\alpha,\beta;\eta)= \sum_{\zeta \in \Es} \chi(\alpha \zeta + \beta
{\zeta}^{-1}) \eta(\zeta).
\end{equation*}
In particular, we write $K(\alpha,\beta)$ for
$K(\alpha,\beta;\eta_1)$, the (standard) Kloosterman sum.

For any $\eta \in \hat{\Es}$, we define the Gauss sum $G(\eta)$
 ($=G_{n,q}(\eta)$) over $E$ by
\[ G(\eta):= \sum\limits_{w \in \Es} \chi(w) \eta(w). \]
 It is clear that some Kloosterman sums will
reduce to Gauss sums.

In what follows, we will use the following properties of
Kloosterman and Gauss sums.  For further details, the reader is
referred to \cite{Coh200} or a reference book such as
\cite{LiNi86}.

\begin{lemma} \label{Kloost bounds}
Let $\eta$ be a multiplicative character of $E$.  Then
\[ K(0,0;\eta)= \left\{ \begin{array}{ll}
q^n-1, & \mbox{if $\eta=\eta_1$,}\\ 0, & \mbox{otherwise.}
\end{array}\right.\]
Further, if either $\eta \neq \eta_1$ or $\alpha$, $\beta \in E$
are not both zero, then
\[ |K(\alpha,\beta; \eta)| \leq 2q^{\frac{n}{2}}. \]
\end{lemma}

\begin{lemma} \label{Kloost}
\begin{itemize}
\item[(i)] If $\alpha$ ($\neq 0$), $\beta \in E$, then
\[ K(\alpha,\beta; \eta)=\bar{\eta}(\alpha)K(1,\alpha
\beta;\eta). \] \item[(ii)] If $\beta \neq 0$, then
$K(0,\beta;\eta)=\eta(\beta)G(\bar{\eta})$. \item[(iii)] If $\alpha
\neq 0$, then $K(\alpha,0;\eta)=\bar{\eta}(\alpha)G(\eta)$.
\end{itemize}
\end{lemma}

\begin{lemma} \label{Gauss}
\begin{itemize}
\item[(i)] $G(\eta_1)=-1$. \item[(ii)] If $\eta \neq \eta_1$,
then $|G(\eta)|=q^{\frac{n}{2}}$.
\end{itemize}
\end{lemma}

\begin{prop} \label{formulation1}
Assume that $m$ is a divisor of $Q$, and $g$, $h$ are divisors of
$x^n-1$. Then
\begin{equation*}
N(m,g,h)=\theta(m) \Theta(g) \Theta(h)  \int\limits_{d|m}
\int\limits_{D_1|g} \int\limits_{D_2|h} K(\delta_{D_1},\delta_{D_2};
\eta_d).
\end{equation*}
\end{prop}
\begin{proof}   Using the characteristic functions derived above, we
have
\begin{equation}\label{eq:N(m,g,h)}
 N(m,g,h) = \sum_{w \in \Es} \left( \theta(m)
\int\limits_{d|m} \eta_d (w)\right) \left( \Theta(g)
\int\limits_{D_1|g} \chi_{\delta_{D_1}}(w)\right) \left( \Theta(h)
\int\limits_{D_2|h} \chi_{\delta_{D_2}}(w^{-1})\right).
\end{equation}
Thus
\[ N(m,g,h)=\theta(m) \Theta(g) \Theta(h)
\int\limits_{d|m}\int\limits_{D_1|g} \int\limits_{D_2|h} \sum_{w
\in \Es} \chi(\delta_{D_1} w + \delta_{D_2} w^{-1}) \eta_d(w), \]
and the result follows from the definition of the generalized
Kloosterman sum.
\end{proof}

>From this, we obtain the following expression.

\begin{prop} \label{formulation2}
Assume that $m$ and $g$ are divisors of $Q$ and $x^n-1$
respectively.  Then
\begin{eqnarray} \label{Neq}
N(m,g,h) &=& \theta(m) \Theta(g) \Theta(h)  \times \nonumber\\
 && \Big( q^n + \epsilon + \int\limits_{\substack{d|m,\\ d \neq 1}}
\int\limits_{\substack{D_1|g,\\ D_1 \neq 1}}
\bar{\eta_d}(\delta_{D_1})G(\eta_d) + \int\limits_{\substack{d|m,\\
d \neq 1}} \int\limits_{\substack{D_2|h,\\ D_2 \neq 1}}
\eta_d(\delta_{D_2})G(\bar{\eta_d}) \nonumber\\ &&+
\int\limits_{\substack{D_1|g,\\ D_1 \neq 1}}
\int\limits_{\substack{D_2|h,\\ D_2 \neq 1}}
K(\delta_{D_1},\delta_{D_2}) + \int\limits_{\substack{d|m,\\ d \neq
1}} \int\limits_{\substack{D_1|g, \\D_1 \neq 1}}
\int\limits_{\substack{D_2|h,\\ D_2 \neq 1}}
K(\delta_{D_1},\delta_{D_2}; \eta_d) \Big),
\end{eqnarray}
where
\[ \epsilon = \left\{ \begin{array}{lll}
-1, & \mbox{if $g=h=1$,}\\ +1, & \mbox{if $g \neq 1$ and $h \neq 1$,}\\
0, & \mbox{otherwise.}
\end{array}\right.\]
\end{prop}
\begin{proof}  We combine the formulation of Proposition
\ref{formulation1} with the results of Lemma \ref{Kloost bounds},
Lemma \ref{Kloost} and Lemma \ref{Gauss}.  If $d=1$, then the
Kloosterman sum takes the value $q^n-1$ when $D_1=D_2=1$,
$\eta_1(\delta_{D_2})G(\bar{\eta_1})=-1$ when $D_1=1$ and $D_2 \neq
1$, and $\bar{\eta_1}(\delta_{D_1})G(\eta_1)=-1$ when $D_2=1$ and
$D_1 \neq 1$.  If $d \neq 1$, then we obtain a contribution of $0$
when $D_1=D_2=1$, $\bar{\eta_d}(\delta_{D_1})G(\eta_d)$ when $D_1
\neq 1$ and $D_2=1$, and $\eta_d(\delta_{D_2})G(\bar{\eta_d})$ when
$D_1=1$ and $D_2 \neq 1$.  Note that the $\epsilon$ term in the
statement of the result arises from the situation when $d=1$,
$D_i=1$ and $D_j \neq 1$ (where $\{i,j\}=\{1,2\}$); for example in
the ``$D_1=1$" case we have $\int\limits_{D_2|h, D_2 \neq 1}
(-1)=-\sum_{D_2|h, D_2 \neq 1} \mu(D_1)$, which takes value $0$ when
$h=1$ and $1$ when $h \neq 1$.
\end{proof}

>From Proposition \ref{formulation2} and the size of the Kloosterman
and Gauss sums, we immediately derive a lower bound for $N(m,g,h)$.
Write $W(m)=2^{\omega(m)}$ for the number of square-free divisors of
$m$, where $\omega$ counts the number of distinct primes in $m$, and
similarly define $W(g)$.

\begin{cor} \label{Nbound}
Under the conditions and with the notation of Proposition
$\ref{formulation2}$,
\begin{eqnarray*}\label{ineq:N(m,g,h)}
 N(m,g,h) & \geq & \theta(m) \Theta(g) \Theta(h) (  q^n +
 \epsilon\\
 &-& q^{n/2}[2W(m)W(g)W(h) - (W(m)+1)(W(g)+W(h)) +2]. )
\end{eqnarray*}
In the case when $g=h$, this inequality takes the form
\begin{equation}\label{ineq:N(m,g,g)}
 N(m,g,g)\geq \theta(m) {\Theta(g)}^2 \left(  q^n + \epsilon_g
- 2q^{n/2}(W(m)W(g)-1)(W(g)-1) \right)
\end{equation}
where
\[ \epsilon_g = \left\{ \begin{array}{ll}
-1, & \mbox{if $g=1$,}\\ +1, & \mbox{if $g \neq 1$.}
\end{array}\right.\]
In particular,
\begin{equation}\label{UsualIneq}
N(m,g,h) \geq \theta(m)\Theta(g)\Theta(h)
q^{n/2}(q^{n/2}-2W(m)W(g)W(h)).
\end{equation}
\end{cor}
\begin{proof} The bounds of Lemmas \ref{Kloost bounds} and
\ref{Kloost} yield for the  sum of the ``integrals" in the identity
($\ref{Neq}$) the absolute bound
\[2(W(m)-1)(W(g)-1)(W(h)-1)+2(W(g)-1)(W(h)-1) +(W(m)-1)(W(g)+W(h)-2).\]
Rearrangement gives the result.
\end{proof}

The following simple bound for $W(m)$, the number of square-free
divisors of $m \in \mathbb{N}$, will be useful in what follows.
The proof is immediate using multiplicativity.
\begin{lemma} \label{Wbound}
For any positive integer $m$,
\begin{equation} \label{eq:Wbound}
W(m)\leq c_m m^{1/4},
\end{equation}
where $c_m= \frac{2^s}{(p_1 \ldots p_s)^{1/4}},$ and
$p_1,\ldots,p_s$ are the distinct primes less than $16$ which
divide $m$.\\ In particular,  for all $m \in \mathbb{N}$,
$c_m<4.9$, and for all odd $m$, $c_m<2.9$.
\end{lemma}

In what follows we recall the notation $\n$ defined by $n=p^b\n,\ p
\nmid \n$.
\begin{prop}\label{small n^*}
Let $q$ be a prime power and let $n(\geq 3) \in \mathbb{N}$ with $\n
\leq 4$. Suppose, in addition,  $q \equiv 2\, (\mod \, 3)$ if
$\n=3$, and $q \equiv 3\, (\mod \, 4)$ if $\n=4$. The pairs
$(q,n)=(2,3),\ (2,4)$ and $(3,4)$ are not PFF. Otherwise, $(q,n)$ is
a PFF pair.
\end{prop}
\begin{proof}
We have $Q(q,n)<q^n/\gcd(n,q-1)$, where,  under the given
conditions,
\[ \gcd(n,q-1)= \left\{ \begin{array}{ll}
1, & \mbox{if $\n=1$ or $3$,}\\ 2, & \mbox{if $\n=2$ or $4$.}
\end{array}\right.\]
Moreover,  $N(Q,x^n-1,x^n-1)=N(Q,g(x),g(x))$, where $g$ factorizes
into $F$-irreducibles as
\[ g(x) = \left \{ \begin{array}{rll}
x-1, & \mbox{if $\n=1$} & \mbox{or $\n=n=3$}, \\ (x-1)(x+1),&
\mbox{if $\n=2$} & \mbox{or $\n=n=4$,} \\ (x-1)(x^2+x+1), &
\mbox{if $\n=3\,<n$,} &\\ (x-1)(x+1)(x^2+1), & \mbox{if
$\n=4\,<n$,}&
\end{array} \right. \]
using Lemma \ref{Red2} when $n=3$ or $4$.  It follows from
Corollary \ref{Nbound} and Lemma \ref{Wbound} that $N:=
N(Q,x^n-1,x^n-1)$ is positive whenever
\begin{equation}\label{ineq:orig small n^*}
q^{\frac{n}{2}}>2(W(Q)W(g)-1)(W(g)-1)-\epsilon_g q^{-\frac{n}{2}},
\end{equation}
and hence whenever
\begin{equation}\label{ineq:small n^*}
\left ( q^n(q-1)\right)^{1/4}>Ac_Q,
\end{equation}
 where \[A= \left \{
\begin{array}{rll} 4, & \mbox{if $\n=1$} & \mbox{or
$\n=n=3$}, \\ 2^\frac{11}{4} \cdot 3,& \mbox{if $\n=2$} & \mbox{or
$\n=n=4$,}
\\ 24, & \mbox{if $\n=3 \, <n$,} &\\ 2^\frac{15}{4}\cdot 7 & \mbox{if
$\n=4 \, <n$.}&
\end{array} \right. \]
We now consider when (\ref{ineq:small n^*}) holds for each of the
values of $A$, using an appropriate bound for $c_Q$.   We use
notation like $(q_0+,n_0+)$ to signify any pair $(q,n)$ with $q
\geq q_0,\; n\geq n_0$.

 \vspace{2mm}
 \noindent $\bullet \quad$  Assume $A=4$.   Then (\ref{ineq:small n^*}) holds
 with $c_Q<4.9$ for $(3+,11+)$, $(4+,8+)$, $(5+,7+)$, $(7+,6+)$, $(8+,5+)$, $(13+,4+)$, $(23+,3+)$; with
 $c_Q < 2.9$ for $(2,15+)$; and with $c_Q < 3.2$ for $(3,9)$ (when $3\nmid Q$). For $n^*=1$, direct
 application of inequality (\ref{ineq:orig small n^*}) establishes the result for
 $(5,5)$, $(8,4)$  and $(4,4)$ (for this last, (\ref{ineq:orig small n^*})
 reduces to $16>14$), leaving only the pairs  $(2,4)$, $(2,8)$ and
 $(3,3)$. When $q=2$, one of the sole reciprocal pair of primitive  quartics has zero trace so there does not exist a PFF polynomial.
 Otherwise, a PFF polynomial for the  case (2,8) is given in
 Section $\ref{two}$; one for (3,3) is in Section $\ref{three}$.
 For the case $n=n^*=3$, inequality (\ref{ineq:orig small n^*}) establishes the result for $(17,3)$,
 $(11,3)$, $(8,3)$ and $(5,3)$. When $q=2$ one of the pair of primitive  cubics has zero trace so there does not exist a PFF polynomial.

 \vspace{2mm}
 \noindent $\bullet \quad$  Assume $A=2^{\frac{11}{4}} \cdot 3$.  Then (\ref{ineq:small n^*}) holds
 with $c_Q<4.9$ for $(3+,17+)$, $(4+,13+)$, $(5+,11+)$, $(6+,10+)$, $(7+,9+)$, $(8+,8+)$, $(11+,7+)$,
 $(14+,6+)$, $(22+,5+)$, $(40+,4+)$; and with $c_Q < 3.2$ for $(3,15+)$.  For the case $n^*=2$, direct
 application of inequality (\ref{ineq:orig small n^*}) establishes the result for $(5,10)$ and $(9,6)$,
 leaving only the pair $(3,6)$.  When $n^*=n=4$, inequality (\ref{ineq:orig small n^*}) establishes the
 result for $(31,4)$, $(27,4)$ and $(23,4)$ (for which $W(Q) \leq 16$) and $(19,4)$ ($W(Q) \leq 8$).
 This leaves pairs $(11,4)$ and $(7,4)$, $(3,4)$.  When $q=3$ there are $4$ primitive quartics with non-zero
traces, namely $f(\pm x)$  where $f(x)=x^4+x^3+x^2-x-1$, together
with their reciprocals. None is a PFF polynomial.  On the other
hand, direct verification yields PFF polynomials as follows.

\bigskip
$
\begin{tabular}{c|c}
  \hline
   $(q,n)$ & PFF polynomial\\ \hline \hline
  $(11,4)$ & $x^4+x^3-5x+2$ \\ \hline
  $(7,4)$ & $x^4+x^3-x^2-x-2$ \\ \hline
\end{tabular}
$

\bigskip
\vspace{2mm}
 \noindent $\bullet \quad$  Assume $A=24$.   Then (\ref{ineq:small n^*}) holds
 with $c_Q<4.9$ for $(16+,6+)$, $(5+,12+)$, and $(2,25+)$.  Inequality (\ref{ineq:orig small n^*})
 establishes the result for $(8,6)$: for $q=2$, degrees $6, 12$ and
 $24$ are treated in Section $\ref{two}$.

\vspace{2mm}
 \noindent $\bullet \quad$  Assume $A=2^{\frac{15}{4}} \cdot 7$.   Then (\ref{ineq:small n^*}) holds
 with $c_Q<4.9$ for $(7+,12+)$, $(4+,20+)$ and $(3+,22+)$.  This leaves the pair
 $(3,12)$ which is treated in Section $\ref{three}$.

\end{proof}

\section{The sieve}\label{Sieve}

In this section, we introduce our key tool: a sieve with both
additive and multiplicative components.  For a given pair $(q,n)$,
let $m|Q$, $f |x^n-1$ and $g|y^n-1$. Let $m_1, \ldots, m_r$ be
factors of $m$, for some $r \geq 1$, and let $f_1, \ldots, f_r$
and  $g_1, \ldots, g_r$  be factors of $f$ and $g$ respectively.
We call $\{(m_1,f_1,g_1) \ldots, (m_r,f_r,g_r)\}$
 a set of \emph{complementary divisor triples} of $(m,f,g)$ with common
divisor triple $(m_0,f_0,g_0)$ if the primes in $\mathrm{lcm}
\{m_1, \ldots, m_r\}$ are precisely those in $m$, the irreducibles
in $\mathrm{lcm} \{f_1, \ldots, f_r\}$ are precisely those in $f$,
the irreducibles in $\mathrm{lcm} \{g_1, \ldots, g_r\}$ are
precisely those in $g$ and, for any distinct pair $(i,j)$, the
primes and irreducibles in $\gcd(m_i,m_j)$, $\gcd(f_i,f_j)$ and
$\gcd(g_i,g_j)$ are precisely those in $m_0$, $f_0$ and $g_0$
respectively.  Observe that the value of $N(m,f,g)$ depends only
on the primes and irreducibles present in $m$, $f$ and $g$. The
following result extends Theorem 3.1 of \cite{Coh100}.

\begin{prop}[Sieving inequality]\label{sieve} For divisors
$m$ of $Q$, $f$ of $x^n-1$ and $g$ of $y^n-1$, let
$\{(m_1,f_1,g_1),\ldots, (m_r,f_r,g_r)\}$ be a set of
complementary divisor triples of $(m,f,g)$ with common divisor
triple $(m_0,f_0,g_0)$. Then
\begin{equation}\label{ineq:sieve}
N(m,f,g) \geq \left (\sum_{i=1}^r N(m_i,f_i,g_i)\right) -
(r-1)N(m_0,f_0,g_0).
\end{equation}
\end{prop}
 \begin{proof}  When $r=1$, the result is trivial. For $r=2$, denote the set of
elements $w \in \Es$ such that $w$ is $m$-free and $f$-free and
$w^{-1}$ is $g$-free, by ${\cal S}_{m,f,g}$. Then ${\cal
S}_{m_1,f_1,g_1} \cup {\cal S}_{m_2,f_2,g_2} \subseteq {\cal
S}_{m_0,f_0,g_0}$, while ${\cal S}_{m_1,f_1,g_1} \cap {\cal
S}_{m_2,f_2,g_2} = {\cal S}_{m,f,g}$, and the inequality holds by
consideration of cardinalities. For $r \geq 2$, use induction on
$r$.
\end{proof}

We observe that, in Proposition \ref{sieve}, $mfg$ can be regarded
as a formal product whose ``atoms" are either prime factors of $Q$
or irreducible factors of $x^n-1$ or $y^n-1$.     Write $k$ for the
(radical of) $mfg$ and $k_0$ for (that of) $m_0f_0g_0$; we shall
refer to $k_0$ as the {\em core} of $k$. Also write $N(k)$ for
$N(m,f,g)$ (so that, in a natural sense, $W(k)=W(m)W(f)W(g)$).
Consider an application of the sieve in which, for each $i=1,
\ldots,r$, $m_if_ig_i$ runs through the values of $k_0p_i$ as $p_i$
runs through atoms of $k$ not in $k_0$. We shall call this a
\emph{$(k_0,r)$ decomposition} of $k$. Given a $(k_0,r)$
decomposition, define $\delta=1-\sum_{i=1}^r\frac{1}{|p_i|}$ with
$|p|=p$ when $p$ is a prime (integer) and $|p| = q^{\deg p}$ when
$p$ is an irreducible polynomial and set  $\Delta =
\frac{r-1}{\delta}+2$.  As we shall see, \emph{it is crucial that
$\delta$ is positive} for the $(k_0,r)$ decomposition selected. In
particular, when $r=1$ (the \emph{non-sieving situation}),  then
($\ref{ineq:sieve}$) is a trivial equality, $W(k)=2W(k_0)$ and
$\Delta =2$.
\begin{prop} \label{keycriterion} In the above notation, for a given pair $(q,n)$, let  $k$ denote the formal product
$mfg$, where $m|Q,\ f|x^n-1$ and $g|y^n-1$.  Suppose that
\begin{equation}\label{eq0}q
> (2W(k))^{2/n}.
\end{equation} Then $N(k)$ is positive.

More generally, for a $(k_0,r)$ decomposition as described above,
suppose that $\delta$ is positive and

\begin{equation}
\label{eq1} q > (2W(k_0) \Delta)^{2/n}.
\end{equation}
Then $N(k)$ is positive.
\end{prop}
\begin{proof} The non-sieving criterion ($\ref{eq0}$) follows
immediately from (\ref{UsualIneq}) of Corollary \ref{Nbound}.

For ($\ref{eq1}$), define $\Theta(k)= \theta(m)\Theta(f)\Theta(g)$
and  write ($\ref{ineq:sieve}$) in the form
\begin{eqnarray} \label{varsieve}
N(k) & \geq & \delta
N(k_0)+\sum_{i=1}^r\left(N(k_0p_i)-\left(1-\frac{1}{p_i}\right)N(k_0)\right)\nonumber\\
 &=& \delta\Theta(k_0)\left(q^n+ \sum_{\substack{d|k_0\\ d \neq 1}}U(d)\right)\
  +\ \Theta(k_0)\sum_{i=1}^r \left(1-\frac{1}{p_i}\right)\sum_{\substack{d|k_0p_i\\ d \nmid
 k_0}}U_i(d),
\end{eqnarray}
where the sums over $d$ are over ``square-free" formal factors of
the formal products $k_0$ and $k_0p_i$ and, by the estimates of
Lemmas \ref{Kloost bounds} and \ref{Gauss} (as already used in
Corollary \ref{Nbound}), each of the expressions $U(d)$ and $U_i(d)$
in absolute value  do not exceed $2q^{n/2}$. Granted that
$\delta>0$, it follows that $N(k)$ is positive whenever
\[ \delta q^{n/2} > 2\delta
W(k_0)+2\sum_{i=1}^r(W(k_0
p_i)-W(k_0))\left(1-\frac{1}{p_i}\right).\] The result follows since
$W(k_0p_i)-W(k_0)=W(k_0)$ and
$\displaystyle{\sum_{i=1}^r\left(1-\frac{1}{p_i}\right)=r-1+\delta}$.
\end{proof}

In applying (\ref{eq1}) to the PFF problem, $k$ is taken to be
$Q(x^n-1)(y^n-1)$; in fact, by the discussion in Section
\ref{sec:2} we may take $k=Q(x^{\n}-1)(y^{\n}-1)$.  Generally, we
take $g_0(x)=f_0(x)$, although if necessary, a more general set of
``complementary divisor triples" or the full form of Corollary
\ref{Nbound} can be used.

We illustrate the direct use of the sieve in dealing with the case
when $n^*=q-1$.

\begin{prop}\label{n=q-1}
Let $q (\geq 4)$ be a prime power and $n (\geq 3) \in \mathbb{N}$.
Suppose $n^*=q-1 >2$.  The pairs $(q,n)=(5,4)$ and  $(4,3)$ are not
PFF. Otherwise, $(q,n)$ is a PFF pair.
\end{prop}

\begin{proof}
We use a $(k_0,r)$ decomposition of $k=Q(x^{\n}-1)(y^{\n}-1)$.
Here $Q=\frac{q^n-1}{(q-1)^2}$ and all polynomial atoms are
linear.

As a first step, we use the additive sieve (alone) with
$f_0(x)=g_0(x)$.  Clearly $\frac{x^{\n}-1}{f_0(x)}$ and
$\frac{y^{\n}-1}{g_0(y)}$ have the same number, $l$ say, of (linear)
factors. To ensure that $\delta$ is positive, of necessity $2l < q$.
Specifically, for $q$ odd (whence $\n$ even) take $l={\n}/2$.  Then
$\delta=1- \frac{\n}{q}=\frac{1}{q}$ and $\Delta=
\frac{\n-1}{\delta} +2={\n}^2+1$. Moreover, $W(f_0)=W(g_0)=
2^{\n/2}$. Thus $(\ref{eq1})$ becomes
\begin{equation}\label{dagger42}
q^{n/2}> 2^{\n+1}({\n}^2+1)W(Q).
\end{equation}

Otherwise, for $q$ even (whence $\n$ odd) take $t=({\n}-1)/2$. Then
$\delta =\frac{2}{q}$ and $\Delta= \frac{\n-2}{\delta}
+2=\frac{(\n-2)(\n+1)}{2}+2 = \frac{{\n}^2-\n+2}{2} <
\frac{{\n}^2+1}{2}$. Now, $W(f_0)=W(g_0)= 2^{(\n+1)/2}$.
Accordingly, $(\ref{dagger42})$ remains a valid sufficient
condition.

 By Lemma $\ref{Wbound}$, $W(Q) < c_Q \frac{q^{n/4}}{\sqrt{q-1}}$. Hence we obtain the sufficient
condition
\begin{equation}\label{star42}
 q^{n/4} > \frac{c_Q 2^q((q-1)^2+1)}{\sqrt{q-1}}.
\end{equation}

First assume that $n=\n$.  Then inequality (\ref{star42}) is
satisfied whenever $n =q-1 \geq 37$. Therefore we can suppose $q
\leq 37$. Next, since $q=n+1 \leq 37$, a straightforward
calculation yields that $\omega(Q) \leq 33$. Now (\ref{dagger42})
yields the sufficient condition
 \[\frac{(n+1)^{n/2}}{2^{n+33}} > 2(n^2+1).\]
This is satisfied whenever $n \geq 26$ ($q \geq 27$).  We may
therefore assume that $q \leq 25$. Another repetition of the
additive sieve (without factorization of $q^n-1$) disposes of
$q=25$.  Next, we introduce a non-trivial multiplicative component
to the sieve (i.e., $m_0 \neq Q$). Factorize $Q$ and take $m_0$ to
be the product of all the primes in $Q$ which are less than $q$
(these are ``worse" than all the linear polynomials $p$ in $k$ since
the latter have $|p|=q$). This deals with $16 \leq q \leq 23$ (or
$25$).  We illustrate in the case when $q=17$. Here $n=16$ and $Q$
has prime factors $3, 5, 29, 18913, 41761, 184417$, so that
$m_0=15$.  Take $t=6$. Then $r=16$, $\delta=
1-\frac{1}{29}-\frac{1}{18913}- \frac{1}{41761}
-\frac{1}{184417}-\frac{12}{17}>0.2595\ldots$, $\Delta=59.79\ldots$
and $W(k_0)=W(m_0)W(g)^2=2^2 \cdot 2^{20}=2^{22}$.  Hence
$(W(k_0)\Delta)^(2/n)<12.24<17$.

Direct verification deals with five of the seven remaining cases ($7
\leq q \leq 13$): see table below. On the other hand, when $q=5$,
given a root $\alpha$ of any of the 32 primitive quartics over
$F=GF(5)$ for which the coefficients of $x^3$ and $x$ are both
non-zero, either $\alpha$ or $1/\alpha$ is \emph{not} free over $F$.
Hence (5,4) is not a PFF pair.  Similarly, when $q=4$, none of the
12 primitive cubics is a PFF polynomial.

In the case when $n>\n$, condition (\ref{star42}) is satisfied for
$q>11$ with  $n \geq 2\n$, for $q > 7$ with $n \geq 3\n$, for $q
>4$ with $n \geq 5\n$, and for $q=4$ (whence $c_Q=2.9$) for $n \geq 8\n=24$.
 The only pairs not covered by this are
$(8,14)$, $(4,12)$ and $(4,6)$.  For $(8,14)$ direct substitution in
condition (\ref{dagger42}) yields the result. For $(4,12)$, use
$(\ref{eq1})$ with multiplicative sieving alone. Specifically,
$Q=5\cdot 7\cdot 13 \cdot 17\cdot 241$. Take the core to be
$(x^3-1)(y^3-1)$ and let all the $t=5$ primes in $Q$ be sieving
primes.  Then $\delta >0.5172$ and $(2W(k_0)\Delta)^{2/n}<3.29<4$.
Finally, a PFF polynomial of degree $6$ is given in Section
$\ref{four}$.

To complete the proof here is the  promised  table of PFF
polynomials.
\bigskip
\bigskip

\begin{tabular}{c|c|c}
  \hline
  $(q,n)$ & PFF polynomial&  polynomial for $u$ \\ \hline \hline
  $(13,12)$ & $x^{12}+x^{11}-3x+2$&\\ \hline
  $(11,10)$ & $x^{10}+x^9-2x+2$ &\\ \hline
  $(9,8)$&\tiny{$x^8-(u-1)x^7-x^6-x^5-(u+1)x^4+(u-1)x^3+(u+1)x^2-x-u$}& $u^2-u-1$\\ \hline
  $(8,7)$& \!\!\tiny{$x^7+x^6+(u+1)x^5+(u^2+1)x^4+(u^2+u+1)x^3+u^2x^2+ux+u^2+u$}\!\!\!\!&$u^3+u+1$\\ \hline
  $(7,6)$&$x^6+x^5+x^2-x+3$&\\ \hline

        \end{tabular}

\end{proof}

\bigskip
\bigskip

\subsection{Key strategy: applying the sieve in the general
case}\label{keys}

In this section, we derive an inequality which provides a sufficient
condition for a pair $(q,n)$ to be a PFF pair in the general case,
by considering a specific factorisation of $x^{n^*}-1$ followed by a
``core-atom" application of the sieve.  The universal value  of this
strategy  can be judged from the fact, in what follows, only a
single case, namely  $(2,21)$, arose for which another factorisation
 succeeded where the key strategy failed.   While the sieve has
both an additive and multiplicative component, we note that it is
often possible to obtain our desired result by using the additive
part alone; correspondingly, we state two versions of our main
inequality. The multiplicative part of the sieve is a useful tool in
dealing with cases where the value of $q$ is small.

Denote by $s$ the positive integer $\mathrm{ord}_{n^*}q$, i.e.
$n^*|q^s-1$ with $s$ minimal; then every irreducible factor of
$x^{n^*}-1$ over $F$ has degree dividing $s$.  Write $x^{n^*}-1$ as
$g(x)G(x)$, where $G$ is the product $\prod_{i=1}^{r}G_i$ of the
($r$, say) irreducible factors ($G_1,\ldots,G_r$, say) of degree
$s$, and $g$ is the product of those with degree less than $s$ (with
$g=1$ if $s=1$). Let $m:=\mathrm{deg}\, g$.  Note that
$r=\frac{n^*-m}{s}$.  For the next result suppose that the set of
$\w(Q)$ distinct prime divisors of $Q$ is partitioned into a set of
$t$ ``sieving" primes $\{l_1, \ldots, l_t\}$ and a set of $u$ primes
whose product is the multiplicative core  $m_0$. Thus $t+u=\w(Q)$;
in particular $u=\w(Q)$ when there is no multiplicative sieving.

\begin{prop} \label{KeyIneq}
Assume the notation defined above.  Then $N(Q,x^n-1,y^n-1)>0$
whenever
\begin{equation}\label{KeyIneqA}
 q^{n/2}>2^{1-t}W(Q)W(g)^2 \left(
\frac{q^s(2(n^*-m)+s(t-1))}{sq^s(1-\sum_{i=1}^t \frac{1}{l_i})
-2(n^*-m)}+2\right),
\end{equation}
provided the displayed denominator in the right side of
$(\ref{KeyIneqA})$ is positive.

In the case of additive sieving only, we have the sufficient
condition
\begin{equation}\label{KeyIneqB}
q^{n/2}>2 W(Q) W(g)^2 \left( \frac{q^s(2(n^*-m)-s)}{sq^s-2(n^*-m)}+2
\right),
\end{equation}
provided the denominator in $(\ref{KeyIneqB})$ is positive.
\end{prop}
\noindent{\bf Note.}  Since $\n|q^s-1$ the denominator in
$(\ref{KeyIneqB})$  is always positive unless $s=1$ and $\n=n=q-1$
(which case is covered by Proposition $\ref{n=q-1}$).

\begin{proof} Take  $2r+t$
complementary divisors with core $k_0=m_0g(x)g(y)$, namely
$\{k_0G_i(x)\ ,i=1,\ldots,r \}$ , $\{ k_0 G_i(y),\ i=1,\ldots,r \}$
and $\{ k_0 l_i, i=1,\ldots,t \}$.

Then $N(Q,x^n-1,y^n-1)$ is positive, by (\ref{eq1}), if
\[ q^{n/2}>2 W(m_0) W(g)^2 \left(
\frac{2r+t-1}{1-\sum_{i=1}^t
\frac{1}{l_i}-\sum_{i=1}^{2r}\frac{1}{q^s}}+2 \right),\] i.e., if
\[ q^{n/2}>2 \cdot 2^u W(g)^2 \left(
\frac{2rs + s(t-1)}{s(1-\sum_{i=1}^t \frac{1}{l_i})-
\frac{2rs}{q^s}}+2 \right),\] i.e., since $rs=n^*-m$, if
$(\ref{KeyIneqA})$ holds.

\end{proof}

\section{Some special cases}

Before treating the problem in its most general setting, we give
separate consideration to some special cases, where the values of
$q$ and $n$ are related, or when $n$ is of a distinguished type
(e.g., prime).

\begin{prop}\label{n|q-1}
Let $q$ $(\geq 5)$ be a prime power and let $n$ $(\geq 3) \in
\mathbb{N}$. Suppose that $\n$ $(>2)$ divides $q-1$ but $\n \neq
q-1$.  Then $(q,n)$ is a PFF pair.
\end{prop}
\begin{proof}
Here we have $G(x)=x^{\n}-1$, $g(x)=1$ and, since $(n,q-1)=\n$, we
have $Q=(q^n-1)/(\n(q-1))$. Moreover $s=1$ and $m=0$.  Note that
here $3 \leq n^* \leq (q-1)/2$; if $n^*<(q-1)/2$, then $n^* \leq
(q-1)/3$.


Inequality (\ref{KeyIneqB}) yields the sufficient condition
 \begin{equation}\label{IneqC}
 q^{n/2} > 2W(Q)\left(\frac{2n^*(q-2)+q}{q-2n^*}\right).
\end{equation}
Using the basic bound $W(Q) < \frac{c_Q q^{(n-1)/4}}{[\n
(1-1/q)]^{1/4}}$ we obtain the sufficient condition
\begin{equation}\label{q>T}
q > \left(2 c_Q
\frac{(2n^*(q-2)+2)}{q-2n^*}\right)^{4/(n+1)}\left(\frac{1}{n^*(1-1/q)}\right)^{1/(n+1)}:=T_1,
\end{equation}
say.  Clearly, $T_1 \rightarrow \infty$ as $\n$ approaches
$\frac{q}{2}$. We shall show that an appropriate upper bound $T_2$
for $T_1$ decreases in the range $3 \leq \n \leq \frac{q-1}{3}$.

Since $q-2\n \geq 1 $ and $n^*(1-1/q) >1$, to begin to analyse
$(\ref{q>T})$, we can replace it by the weaker sufficient condition

\begin{equation}\label{51eq1}
q>(2c_Q q(2\n+1))^{4/(n+1)}:=T_2,
\end{equation}
say.

We first consider the case when $n=\n$.  We begin by assuming that
$n \geq 10$: thus  $q \geq 23$.  Taking natural logarithms,
\[ \log T_2= \frac{4}{n+1}( \log(2c_Q q)+
\log(2n+1)).\] For fixed $q$, differentiating with respect to $n$ we
obtain
\[ \frac{d}{dn} \log T_2=-\frac{4}{(n+1)^2}\left(\log
2c_Qq(2n+1)-\left(1+\frac{1}{2n+1}\right)\right),\] which is
negative since  $\log (4n+2)  > 1 + \frac{1}{2n+1}$ for all $n\geq
1$. So, in the range $10 \leq n \leq \frac{q-1}{2}$, the maximal
value of $T_2$ is attained at $n=10$:  it is certainly less than $q$
for $q\geq 23$.

Now assume $3 \leq n \leq 9$.  Since $q-2n \geq q-18$, we can
replace   (\ref{q>T}) by
\[q>\left(\frac{2c_Q q(2\n+1)}{q-18}\right)^{4/(n+1)}:=T_3,\]
say.  Taking logarithms and differentiating, we find that $T_3$ is a
decreasing function if
\[\log \frac{14 c_Q q}{q-18}>\frac{8}{7},\]
which holds for $q>18$ (since $\log 14>8/7$).  The maximum value of
$T_3$ occurs when $n=3$; it is less than $q$ for $q>14 c_Q+18$,
i.e., $q>86$. This establishes the result except when $q<87$ and $3
\leq n \leq \mathrm{min}(9,\frac{q-1}{2})$.

Using $(\ref{IneqC})$, with the $c_Q$ bound, we find from a
computational check that the result holds for all remaining $(q,n)$
except $(19,9), (17,8), (19,6), 13,6), (16,5), (11,5)$ and
appropriate values of $(q,4),\ n \leq 29$  (5 values) and $(q,3), n
\leq 49$ (9 values). For all remaining values, $\w(Q)\leq 4$; taking
exact  values deals (via  $(\ref{IneqC})$) with all pairs except
$(7,3)$, $(16,3)$, $(9,4)$, $(13,4)$, $(11,5)$, $(13,6)$. Invoking
the multiplicative part of the sieve also, i.e., using inequality
(\ref{KeyIneqA}), yields the results for $(13,6)$ ($Q=7 \cdot 61
\cdot 157$, $m_0=7$) and $(16,3)$ ($Q=7 \cdot 13$, $m_0=7$). Direct
verification establishes the other four cases (see table below).

Now suppose $n>\n$, and replace $n+1$ by $2\n+1$ in (\ref{51eq1}) to
obtain the sufficient condition
\begin{equation}\label{51eq2}
q>(2 c_Q q(2\n+1))^{4/(2\n+1)}:= T_4,
\end{equation}
say. We begin by assuming that
$\n \geq 5$ and $q>13$. Taking logarithms and differentiating,
\[ \frac{d}{d\n} \log T_4=\frac{8}{(2\n+1)^2}(1-\log 2c_Q q
(2\n+1)),\] clearly negative.    So, in the range $5 \leq \n \leq
\frac{q-1}{2}$, the maximum value of $T_4$ is attained at $\n=5$,
and this is less than $q$ for $q>13$.   When $n^*=\frac{q-1}{2}$, we
note that $n \geq 3\n$; using this in condition (\ref{51eq2}), we
find the result holds for $q \geq 9$ (and so in every case).

Finally we consider $3 \leq \n \leq 4$. Since $\n \leq 4$, we can
use a final sufficient criterion, namely
\[q>\left(\frac{q(2\n+1)}{q-8}\right)^{4/(n+1)}:=T_5,\] say. Again
by differentiation, we can check that $T_5$ is a decreasing function
when $q>8$.  The maximum value of $T_5$ occurs when $\n=3$; this is
less than $q$ when $q>13$. This leaves only $n>\n$ with $q=13$,
$\n=3,4$. Using $n \geq 13 \n$ in the sufficient condition yields
the result.
\bigskip

$
\begin{tabular}{c|c|c}
  \hline
  $(q,n)$ & PFF polynomial&  polynomial for $u$ \\ \hline \hline
  $(13,4)$ & $x^4+x^3-x-2$& \\ \hline
  $(11,5)$& $x^5+x^4+3x-2$& \\ \hline
  $(9,4)$&$x^4-x^3+x^2+x-u+1$&$u^2-u-1$ \\ \hline
  $(7,3)$& $x^3+x^2+2x-3$&\\ \hline
\end{tabular}
$

\end{proof}

\bigskip
The following simple lemma  improves  Lemma $\ref{Wbound}$ under the
stated conditions.
\begin{lemma}\label{c=1}
Let $n \geq 5$ be prime, and let $h \in \mathbb{N}$ be squarefree
with each prime divisor of $h$ congruent to $1$ modulo $2n$. Then
\[ W(h)< h^{1/4}, \]
except when $n =5$ and $h=11$.
\end{lemma}

\begin{prop}\label{spec1}
Let $q$ $(\geq 5)$ be a prime power and let $n \in \mathbb{N}$.
Suppose $\n\ (\geq 5)$ does not divide  $q-1$ and either $\n$ is
prime or $\n=q+1$ with $q$ even. Then $(q,n)$ is a PFF pair.
\end{prop}
\begin{proof}
In this case, $x^{n^*}-1$ factors as $(x-1)G(x)$ where $G$ is a
product of $\frac{\n-1}{s}$ factors of degree $s$. We have $s \geq
2$ ($s=2$ if $\n=q+1$); $m=1$, $(n,q-1)=1$ and $Q=\frac{q^n-1}{q-1}$
odd.

By inequality (\ref{KeyIneqB}) of Proposition \ref{KeyIneq}, we have
the sufficient condition (for $N(Q,x^n-1,y^n-1)>0$)
\begin{equation}\label{suff5}
 q^{n/2}-8 W(Q) \left( \frac{2(\n-1)}{s-\frac{2(\n-1)}{q^s}}+1
\right)>0; \end{equation} this certainly holds   if
\[\Delta=\Delta(q,n,s):= (q^n(q-1))^{1/4}- 8c_Q
\left(\frac{2(\n-1)}{s-\frac{2(\n-1)}{q^s}}+1 \right)>0.\]

Concentrating on the ``worst-case scenario" when $n=\n$, we
require
\begin{equation}\label{delta1}
\Delta(q,\n,s)= (q^{\n}(q-1))^{1/4}- 8c_Q
\left(\frac{2(\n-1)}{s-\frac{2(\n-1)}{q^s}}+1 \right)>0.
\end{equation}
In $(\ref{delta1})$  we can take $c_Q<2.9$ since $Q$ is odd. In
fact, when $q$ and $n$ are odd and $n$ is an odd prime, Lemma
\ref{c=1} applies and we can take $c_Q=1$.

Evidently,    $\Delta(q,\n,s)$ is an increasing function of $q$
(with $\n,s$ fixed) and of $s$ (with $q, \n$ fixed).  It is also
increasing with respect to $\n$ with some qualification as regards
to small values of $q, s$.   In fact, with $c_Q=1$, by
differentiation, for given odd $q$ and  $s=2$, $\Delta$ is an
increasing function of $\n$ in the range $(q,\n)=(5+,9+),
(7+,6+),(9+, 5+)$. For even $q$ (take $c_Q=2.9$), the corresponding
pairs are $(8+,6+), (16+, 5+)$.
 For $s=3$,  the pairs need to be $(5+,6+),(7+,5+)$, $q$ odd; $(8+, 6+), (16+, 5+)$, $q$ even.
 For $s\geq4$,   any pair $(5+,5+)$ ($q$ odd) or $(8+,5+)$ is in a region of increasing
$\Delta$.
 Within the above framework, it suffices to establish the result for smallest $q$ and
$n$. It also suffices to take least $s$, i.e., $s=2$.

In the general case, by computation, the result holds for
$(25+,5+)$, $(16+,7+)$, $(9+,9+)$, $(7+,11+)$ and $(5+,17+)$: in
each case within the range of  increasing $\Delta$ with $\n$.

 Suppose first
that $n=\n$.  For the pairs $(q,n)$ not covered by the above, a
number are simply excluded by Lemma \ref{EarlyEx}.  For all but two
remaining pairs, $\Delta(q,n,s)$ is quickly calculated to be
positive; specifically, when
($q,n,s)=(19,5,2),(13,7,2),(11,7,3),(9,7,3),$
$(5,17,5),(5,13,4),(5,17,16), (5,19,9)$ or $(5,23,22)$. The final
two pairs are $(9,5)$ and  $(8,9)$: in each case $s=2$. For these,
$W(Q) =4, 8$, respectively and the sufficient condition
$(\ref{suff5})$ holds.

  Finally, suppose  $n>\n$.  In the $\Delta$ definition and condition (\ref{delta1}),
 replace in the first term $q^{\n}(q-1)$ by $q^{3\n}(q-1)$  ($q$ odd) and by  $q^{2\n}(q-1)$
($q$ even).  Also, set $c_Q =1$ or $2.9$ according as $q$ is odd or
even.  Then, easily, $\Delta(5, \n,2)$ and $\Delta(8,\n,2)$ are
increasing and positive in the respective cases.  This completes the
proof.
\end{proof}

\begin{prop}\label{spec2}
Let $q$ $(\geq 5)$ be an odd prime power and let $n \in \mathbb{N}$.
Suppose $n^*=2l \geq 6$, where either  $l$ is a prime not dividing
$q-1$ or $l=\frac{1}{2}(q+1)$ with  $q \equiv \, 3 \, (\mathrm{mod}
\, 4)$. Then $(q,n)$ is a PFF pair.
\end{prop}
\begin{proof}  When  $l$ is prime then $2 \leq s|l-1$ (since $q^s\equiv
1(\mathrm{mod}\ l)$), whence $\n-2=2(l-1)$ is divisible by $s$.  The
same conclusion holds when $l=\frac{1}{2}(q+1)$, in which case
$s=2$.  Indeed,   in both cases, $x^{\n}-1$ factors into two linear
factors and $\frac{\n-2}{s}$ factors of degree $s$. (Note that
$(\n,q-1)=2$.) Let $\gamma_s=1$ if $s$ is even, or  $2$ if $s$ is
odd: thus, since $l$ divides $\frac{\gamma_sq^s}{2(q-1)}$ then
$l<\frac{\gamma_sq^s}{2(q-1)}$. Apply Proposition  $\ref{KeyIneq}$
with  $m=2$. By inequality (\ref{KeyIneqB}), we have the sufficient
condition
\begin{equation}\label{suff6}
 q^{n/2}-64 W(Q)\left( \frac{\n-2}{s-\frac{(\n-2)}{q^s}}+1
\right) >0.
\end{equation}
 which, as before, is certainly implied by
\[\Delta(q,n,s):= (q^n(q-1))^{1/4}-\frac{64}{2^{1/4}}c_Q \left(
\frac{\n-2}{s-\frac{2\gamma_s}{q-1}}+1 \right) >0.\] Concentrating
on the ``worst-case scenario" when $n=\n$, we require
\begin{equation} \label{delta2}
\Delta(q,\n,s)>0.
\end{equation}
As in Proposition \ref{spec1}, it suffices to establish the result
for smallest $q$ and $n$. We take $s=2$, $\gamma_s=2$ and
$c_Q<4.9$.

By computation, the result holds for $(47+,6+)$, $(23+,8+)$,
$(16+,10+)$, $(11+,12+)$, $(9+,14+)$, $(7+,16+)$ and $(5+,21+)$. We
may now assume that $q \leq 43$.

Suppose first that $n=\n$.  Note that, for $n\geq 14$, the only case
which remains is $(5,14)$. When $n=6$, we find that $W(Q) \leq 2^5$
for all $q<47$ with $q \not \equiv 1(\mathrm{mod} \ 6)$. Using this,
$(\ref{suff6})$  gives the result for $q \geq 19$.   Indeed, for
$q<19$, all except $q=11$ have $W(Q)\leq 2^4$, which gives the
result for $q=17$. This leaves just $q \leq 13$ when $n=6$; in fact,
only $(5,6)$ ($Q=2 \cdot 3^2 \cdot 7 \cdot 19 \cdot 37$)  and
$(11,6)$ ($Q=3^2 \cdot 7 \cdot 31$). Using both the additive and
multiplicative power of the sieve, i.e., using inequality
(\ref{KeyIneqA}), gives the sufficient condition
\[ q^{n/2}>2^{1-t}16W(Q) \left(
\frac{n^*+t-3}{(1-\sum\frac{1}{l_i})-\frac{n^*-2}{q^2}}+2 \right).
\] With $t=3$, this yields the result for $q=11$ ($l_1=7$, $l_2=19$ and
$l_3=37$).  This leaves just $q=5$. When $n=8$, using the
additive-only estimate with $W(Q)=2^3$ and $\gamma_s=1$ gives the
result for $(7,8)$. When $n=10$, all valid $q<16$ have $W(Q)=2^4$;
using this value in the additive-only inequality yields the result
for all $q \geq 7$. Finally, using $W(Q) \leq 2^5$ deals with
$(5,14)$.  Direct verification deals with the remaining case: the
pair $(5,6)$ has PFF polynomial $x^6+x^5+x^3+x^2-x-2$. When $n>\n$,
taking $3\n$ in place of $\n$ in the first term of condition
(\ref{delta2}) yields the result for all pairs.
\end{proof}

\section{Larger fields and degrees}
It is necessary to deal individually with fields of smallest
cardinality, namely $ 2, 3$ and $4$, and their treatment is
deferred to Section $\ref{vsmall}$.  Here we suppose $q \geq 5$.
Even so,  it turns out that $\mathbb{F}_5$ and $\mathbb{F}_7$
require closer attention. From what has been accomplished so far
we may also assume  that $n^* \geq 8$.

We make the following definitions. For $g$ as defined in Section
${\ref{keys}}$, $\omega=\omega(q,n)$ is the number of distinct
irreducible factors of $g$ (so $W(g)=2^{\omega}$), and
$\rho=\rho(q,n)=\frac{\omega(q,n)}{n}$.  For later use, given $n$
also  define $\rho^*=\rho(q,n^*)$, so that $\rho^*/\rho=n/n^*$ is
the power of $p$ in $n$.   As in Section ${\ref{keys}}$, $s$ denotes
the degree of the irreducible factors of $G$. We can suppose that
$s\geq 2$. Also set $n_1:=\gcd(n,q-1)$.

\begin{lemma}[\cite{CoHu03}]\label{rho}
Assume that $n>4$ with  $p \nmid n$ and $q>4$.  Then the following
hold.
\begin{itemize}
\item[\emph{(i)}] If $n=2n_1$ with $q$ odd, then $s=2$ and $\rho=1/2$;

\item[ \emph{(ii)}] if  $n=4n_1$ with $q \equiv 1 (\mathrm{mod} \ 4$), then $s=4$
and $\rho=3/8$;

\item[\emph{(iii)}] if  $n=6n_1$ with $q \equiv 1(\mathrm{mod} \ 6)$, then $s=6$
and $\rho=13/36$;

\item[\emph{(iv)}] otherwise, $\rho \leq 1/3$.
\end{itemize}
\end{lemma}
Because the bounds of Lemma $\ref{rho}$ (taken from \cite{CoHu03})
are  insufficient in themselves when $q=5$ or $7$, there is some
difficulty for these field cardinalities.
  We overcome the obstacle  by a numerical result
related to Lemma $\ref{Wbound}$;  bounds of similar type (such as
Lemma \ref{2W}) will occur in Section $\ref{vsmall}$).

\begin{lemma}\label{1/6bound} Suppose $\omega(h)\geq 49$.  Then
\[W(h) < h^{1/6}.\]

\end{lemma}
\begin{proof} By calculation the result holds when $\omega(h)=49$,
since then $h$ is at least the product of the first $49$  primes.
The result follows since the $50$th prime is $229>2^6$.
\end{proof}

\bigskip
Write the radical of $Q$ as $m_0p_1\ldots p_t$, where $m_0$ is the
core and $p_1, \ldots, p_t$ are the (multiplicative) sieving primes.
When $t=0$ there is no multiplicative sieving. Set $u:=\omega(m_0)$;
thus, often $u= \omega(Q)$.  In  this context, the basic form of
 $(\ref{KeyIneqA})$ in Proposition $\ref{KeyIneq}$ takes the  shape $(\ref{ineq1})$ with
  $(\ref{Rform})$ or $(\ref{Rformadd})$ below (because
$\n-m=\n-\rho n\leq (1-\rho)n$): by contrast,  the \emph{refined
form} does not employ this simplification.

\begin{prop} \label{rhoform}
Suppose that
\begin{equation}\label{ineq1}  q \ >\ R(n),
\end{equation}
 where
\begin{equation} \label{Rform}
R(n)\  =R(n;q)\ = \ \left\{2^{2\rho
n+u+1}\left(\frac{\frac{2(1-\rho)n}{s}
+t-1}{{\delta-\frac{2(1-\rho)n}{sq^s}}}+2\right)\right\}^{2/n}.
\end{equation}
and $\delta= 1-\sum_{i=1}^t\frac{1}{p_i}$ (with $\delta=1$ when
$t=0$). Then $(q,n)$ is a PFF pair.

In particular, when additive sieving alone is being used (i.e.,
$t=0$), then  $R(n)$ takes the  form
\begin{equation}\label{Rformadd}
R(n)\  =R(n;q)\ = \ \left\{2^{2\rho
n+u+1}\left(\frac{\frac{2(1-\rho)n}{s}
-1}{{1-\frac{2(1-\rho)n}{sq^s}}}+2\right)\right\}^{2/n}.
\end{equation}

In the refined form of Lemma $\ref{rhoform}$ both occurrences of
$(1-\rho)n$ are replaced by $\n-\rho n$ in each of $(\ref{Rform})$
and $(\ref{Rformadd})$.

\end{prop}

  Note also that $R(n;q)$
 depends on $q$ (as well as  $n$).   Inasmuch as it is obviously a decreasing function of $q$
 (for fixed values of the other parameters),
 we shall  apply it either  when $q$ has a
specified value or  when $q \geq q_0$ with $q_0$ specified. In what
follows we shall, for convenience of calculation, use alternative
weaker (i.e., larger) forms of $R(n)$ (to be denoted by $R_1(n), \
R_2(n)$, etc): it will be sufficient to show that $(\ref{ineq1})$
holds for the relevant form.

We divide the discussion into two categories according as to whether
$\rho >1/3$ or $\rho \leq 1/3$ as described in Lemma $\ref{rho}$.
When $n^*<n$ then $\rho(q,n) \leq \frac{\rho(q,n^*)}{p} \leq
\frac{\rho(q,n^*)}{2}$.  This means that such pairs fall in the
second category: moreover, from  the size of $\rho(q,n)$,
these scarcely feature  in the discussion.

\begin{prop}\label{bigrho}
Suppose $q \geq 5$ and $n \geq 8$ with $n \nmid(q-1)$. Suppose also
that $\rho(q,n) > 1/3$. Then $(q,n)$ is a PFF pair.
\end{prop}

\begin{proof} The circumstances where $\rho >1/3$ are delineated
in Lemma $\ref{rho}$.   In these, put $n=dn_1$ where  $d=2,4$ or
$6$. Then $Q=\frac{d(q^n-1)}{n(q-1)}$ and $n^* = n < qd$.  By means
of the simple bound $(\ref{eq:Wbound})$  for $W(Q)$ and without
multiplicative sieving, we obtain (as an alternative to $R(n)$)
\begin{equation}\label{bigrhoeq}
R_1(n):\ =\ \left\{c2^{2\rho
n+1}\Big(\frac{d}{n(q_0-1)}\Big)^{1/4}\left(\frac{\frac{2(1-\rho)n}{s}
-1}{{1-\frac{2d(1-\rho)}{sq_0^{s-1}}}} +2\right)\right\}^{4/n}
\end{equation}
(with $c<4.9$ and  $q \geq q_0$)  for  use in $(\ref{ineq1})$.

Because $n^{1/n}$ decreases as $n$ increases, it is seen (with a
little effort)  that $R_1(n)$ decreases as $n\geq 8$ increases under
the given conditions.

>From Lemma $\ref{rho}$, suppose first that $\rho=1/2$ (with $s=2$
and $d=2$ ).  Then $R_1(8;59) < 57$.  Hence $(q,n)$ is a PFF pair
whenever $q \geq 59$.  Indeed, $R_1(12;43)<41.6$,
 and  $R_1(16;37) < 34.7$, etc., thus reducing further
 the list of possible exceptional pairs.  Since $n <2q$, it can
thus  be quickly  checked (using $R_1 $ for $R$ in $(\ref{ineq1})$)
that
 the only pairs not shown to be PFF pairs are $(5,8),\ $  $(7,12),\ $
 $(9,16), \ $  $(11,20),\ $  $(13,8),\ $  $(13,24),\ $  $(17,32), \ $  $(19,12), \ $
 $(19,36), \ $  $(25,16), \ $  $(29,8), \ $ $(31,12), \ $ $(37,8), \ $ $(53,8)$.

These   $14$ pairs were then tested using $(\ref{Rformadd})$,
having calculated $u$ by factorizing $Q$. This was successful
except for $(5,8),\ $ $(7,12),\ $
 $(9,16), \ $  $(13,8)$.  The final stage for these pairs was to
 sieve multiplicatively, also.  Thus, for $(9,16)$, $Q=2\cdot 5
 \cdot 17 \cdot 41 \cdot 193 \cdot 21523361$, the
 largest four primes being the sieving ones.  With $u=2$ this
 yields $R(16) < 7.4$ and hence a PFF pair. Similarly, for
 $(13,8)$, $Q=2\cdot 5
 \cdot 7 \cdot 17 \cdot 14281$, and, again with four sieving
 primes, this yields $R(8) < 11$ and another PFF pair.  This
 process fails, however,  for two  pairs $(5,8)$ and
 $(7,12)$.  For these  we list an explicit PFF polynomial
 as follows.

 \bigskip
 $
\begin{tabular}{c|c}
  \hline
   $(q,n)$ & PFF polynomial \\ \hline \hline
    $(7,12)$ &  $x^{12}+x^{11}-3x-2$\\ \hline
  $(5,8)$ & $x^8+x^7-x^2-x-2 $\\ \hline
\end{tabular}
$

 \bigskip
 Next, suppose from Lemma $\ref{rho}$,  that $\rho=3/8$ (with $s=4$
and $d=4$).  This implies that $n \geq 16$.  We calculate
$R_1(16;19)<17$ and $R_1(13;13)<13$. This excludes only the pairs
$(5,16)$, $(9,32)$ and $(13,16)$. In all these cases, $\w(Q) \leq
7$. Using this in $(\ref{Rformadd})$ with $u=7$, we see that
$(13,16)$ and $(9,32)$ are  (comfortably) PFF pairs. For $(5,16)$,
 use multiplicative sieving. Here $Q=2^2\cdot 3 \cdot 13 \cdot 17
\cdot 313 \cdot 11489$ and we take $u=2, \ t=4$ to yield $\delta=
0.8610$ and $R(16;5)<5$.

\medskip
Finally, suppose from Lemma $\ref{rho}$,  that $\rho=13/36$ (with
$s=6$ and $d=6$). This implies that $n \geq 36$ and  $R_1(36;11))
< 10.9$.   This does leave the pair $(7,36)$ but an application of
$(\ref{Rformadd})$ with $u=11$ yields $R(36;7) < 5$.
 \end{proof}

\bigskip
For the remainder of this section we assume $\rho \leq 1/3$.
Consider the function $R(n;q)$ defined by $(\ref{Rformadd})$. In the
situation to which it   applies,  $s$ and $\rho$ are determined by
$q$ and  $n$.  Nevertheless it is useful sometimes to consider
$R(n;q)$ (and similar expressions) as functions of $n, q, s$ and
$\rho$, more loosely related.  (For instance, since $s\geq 2$ is the
least integer for which $\n$ divides $q^s-1$, then $\n < q^s$ and $s
\leq \phi(\n) < \n$.) It is important to ensure that $sq^s <
2(1-\rho)n$ so that the right side of (\ref{Rformadd}) is a
well-defined positive quantity. It  is a consequence of the next
lemma that, for given $n,q, s$ with $2 \leq s <n$ and $8\leq n<q^s$
(indeed $n< q^2/2$ when $s=2$), then $R(n;q)$ is an increasing
function of $\rho$ for $0 \leq \rho \leq 1/3$.
\begin{lemma}\label{rhoincr} For fixed positive integers $n,q, s$ with $2 \leq s <n$ and
$8\leq n<q^s$ (indeed with $n<q^2/2$ when $s=2$), set
\[ \tau(\rho)=2^{2\rho n}\left(\frac{\frac{2(1-\rho)n}{s}
-1}{{1-\frac{2(1-\rho)n}{sq^{s}}}}\right).\] Then,  $\tau(\rho)$ is
an increasing function for $0\leq \rho \leq 1/3$.
\end{lemma}
\begin{proof} Differentiate to obtain
\begin{equation}\label{diffT}
\tau'(\rho)= K\cdot[\log 2(2(1-\rho)n-s)(sq^s-2(1-\rho)n)\ - \
s(q^s-1)],\
\end{equation}
where $K=\displaystyle{ \frac{ns^2q^{2s}}{(sq^s-2(1-\rho)n)^2}}$ is
a positive function (of all the variables).

If $s=2$ then, since $0 \leq \rho \leq 1/3$ and $n<q^2/2$,
\[\tau'(\rho) \geq \log 2 (\frac{4n}{3}-2)(2q^2-q^2)\ -\ 2q^2
=q^2( \frac{4n}{3}\log2-2)\ -\ 2) \ > \ 0, \]  since $n \geq 8$.

If $3\leq s < n/2$, then, by $(\ref{diffT})$,  for $0\leq \rho \leq
1/3$,
\begin{eqnarray*}
\tau'(\rho)/K &\ \geq\ & n\log 2(\frac{4}{3}-\frac{1}{2})(s-2)q^s-sq^s \\
& \ =\ & q^s\left(s(\frac{5n}{6}\log 2-1)- \frac{5n}{3} \log 2\right)\\
& \ \geq \ &   q^s\left(3(\frac{5n}{6}\log 2-1)- \frac{5n}{3} \log
2\right) \ = \ q^s\left(\frac{5n}{6}\log 2 -3\right)\ > 0,
\end{eqnarray*}
since $n \geq 8$.

Finally, if $ n/2 \leq s \leq n\  (< q^s)$, then, again by
$(\ref{diffT})$,
\begin{eqnarray*}
\tau'(\rho)/K &\ \geq\ &
\frac{n}{3} (sq^s-4s)\log 2)-sq^s= s\left\{\frac{n}{3} \log 2-1)q^s- \frac{4n}{3} \log 2\right\}\\
&\ > \ & s\left[n(\frac{n}{3} \log 2-1)- \frac{4n}{3} \log 2
\right]= \frac{ns}{3}[(n\log 2-3)- 4 \log 2] \ > \ 0,
\end{eqnarray*}
again since $n \geq 8$.
\end{proof}

In practice, it is  convenient to employ a larger  ``starter"
function $\bar{R}(n;q)$, derived from $R(n)$ by taking $\rho=1/3$,
and then using the facts that $n<q^s$ and $s \geq 2$.
\begin{equation} \label{starterR}
(R(n;q)<)\ \bar{R}(n)  =\bar{R}(n;q): = \ \{2^{(2/3)
n+u+1}(2n-1)\}^{2/n}.
\end{equation}
In the result which follows we employ suitable modifications of
these ideas.

\begin{prop}\label{weerho}

Suppose $q \geq 5$ and $n^* \geq 8$ with $n^* \nmid(q-1)$. Suppose
also that $\rho(q,n) \leq  1/3$. Then $(q,n)$ is a PFF pair.
\end{prop}
\begin{proof}
As usual, we generally suppose for simplicity that  $n=\n$ in the
theoretical discussion. Nevertheless, in the treatment of residual
(more delicate) cases, pairs $(q,n)$ for which $n>\n$ are also
considered where relevant.

\bigskip
\noindent{\bf Case O:} $\n=q^2-1$.

In this situation, the argument about $R(n)$ increasing with $\rho$
(to be used elsewhere) fails.  Here $\rho=1/(q+1)$ and $R_1(q^2-1)$
(defined by $(\ref{bigrhoeq})$) has the form
\[ R_1(q^2-1)= (c2^{2q-1}(q^3-q^2-q+2))^{4/(q^2-1)}.\]
With $c=4.9$, it is quickly seen that $R_1(q^2-1)$ decreases and is
less than $9.8$ for $q \geq 11$.  Moreover, when $q=9$, we can take
$c=3.2$ and $R_1(9^2-1)<7.6$ and when $q=8$, we can take $c=2.9$ and
$R_1(8^2-1)<7.8$.  For the pair $(7,48)$, with $s=2$ and $u=13$, we
have $R(48;7) <2.69 <7$.  The discussion of the final pair $(5,24)$
is incorporated with the figures for the most delicate cases in Case
II below. In what follows we  assume (as we may) $\n < q^2/2$ when
$s=2$.

 \noindent {\bf Case I:}  $q \geq 8$.

\noindent Replace $\rho$ by $1/3$ and use Lemma $\ref{Wbound}$ in
$(\ref{Rformadd})$. It therefore suffices that $q > R_2(n)$, where
\begin{equation}\label{weerhoeq}
R_2(n)=R_2(n;q,s)\ =\ \left\{c2^{(2/3)
n+1}\frac{1}{(q-1)^{1/4}}\left(\frac{\frac{4n}{3s}
-1}{{1-\frac{4n}{3sq^{s}}}} +2\right)\right\}^{4/n},
\end{equation}
where $c<4.9$. Here a suitable starter form, derived from
$(\ref{weerhoeq})$ by using  $s \geq 2$ and $n < q^s$ is

\begin{equation}\label{preweerho}
\bar{R}_2(n)= \bar{R}_2(n;q)\ =\ \left\{c2^{(2/3) n+1}\frac{(2n-1)
}{(q-1)^{1/4}}\right\}^{4/n}.
\end{equation}
Evidently $\bar{R}_2(n; q,1/3) )
$ increases as $n$ or $q$ decreases.
With $c=4.9$ , we have $\bar{R}_2(8;49) < 47.5$. Hence the result
holds for $q \geq 49$.

We treat prime powers $q\leq 47$  first by $\bar{R}_2(n)$, to
establish the result for (potentially) large values of $n$ and $s$,
and then by  $R_2(n)$ for more critical values of $n$, with $s$
(close to) $2$. Indeed, to begin, suppose $37 \leq q \leq 47$.  Take
$c=4.9$. Since  $\bar{R}_2(10; 37)< 36$ the result holds for this
range of $q$,
 provided $n \geq 10$.  But also  $R_2(8;37,2)< 32.1$. Hence
the result holds unconditionally.

Smaller values of $q$ are dealt with individually.   For example,
take $q=11$ (so that $c=4.5$ will do).  Then $\bar{R}_2(45;
11)<10.94$, so that we can assume $n \leq 44$ with $n \neq$ a prime
or twice a prime or $12$ (by Propositions $\ref{spec1}$ and
$\ref{spec2}$).  Further, $R_2(37;11,2)< 10.996$, and indeed
$R_2(35;11,3)< 10.6$ (rules out $n=35$), and $R_2(26;11,6)<10.94$
(rules out $n=28, 36$).  We conclude that $n \in
\{8,9,15,16,18,20,21,24,25,30\}$.  For these remaining values,
calculate $u:= \w(Q)$ for use in Proposition \ref{rhoform} by means
of $R(n;11)$ given by $(\ref{Rformadd})$  with $s=2$ and $\rho
=1/3$. In fact, for this set of values of $n$, we have  $u \leq 11$
(attained when $n=24$); indeed, for $n=8, 9$, we have $u \leq 4$.
Now, with $u=11$, we obtain $R(15;11) < 10.7$ and, with $u=4$, we
obtain $R(8; 11)< 9.6$. So the discussion of the case $q=11$ is
therefore complete.

Suppose, next $q=9$ (so that one can take $c=3.2$). Note that we
need also to consider values of $n>\n $ but that, by previous
results and since $\n \geq 8$, it can be supposed that $n \geq 16$.
Since $\bar{R}_2(73; 9) <8.98$, it can be assumed that $n \leq 73$.
Some smaller values of $n$ can be ruled out by  $R_2(n)$. For
example $R_2(56; 9, 3) < 8.9; R_2(55; 9,10)< 8.2; R_2(64;9, 6)< 8.2$
(rules out $n=64,\ 68$); $R_2(49;9,21)< 8$.  The values of $n$ which
remain lie in the set $\{16,20,24,25,28,32,35,36,40,44,48,52,60\}$.
By calculation, $u \leq 17$  (attained at $n=60$); indeed, $u \leq
10$ (attained at $n=40$) except for $n \in \{24, 36, 44, 48, 60\}$.
Finally, take $\rho =1/3, s=2$ in $(\ref{Rformadd})$.  With $u=17$,
we have $R(25; 9) < 8.8$, with $u=11$, then $R(18)<8.6$,  and, with
$u=6$, then $R(16; 9)<6.4$. So the discussion when $q=9$ is
complete.

Finally, suppose $q=8$ (so that one can take $c=2.9$).  The most
delicate degree ($n=9$) has been dealt with in Proposition
\ref{spec1}; more generally, previous results ensure we may assume
$n \geq 15$.  Since $\bar{R}_2(117; 8) < 7.991$ we can suppose that
$n \leq 116$. For $88 \leq n \leq 116$ then $ s \geq 3$ and
$R_2(88;8) < 7.99$ and the result holds. Assume $n \leq 87$, Now
take $\rho=1/3$ and  $s=2$ in $(\ref{Rformadd})$. If $33 \leq n \leq
87$, then $u \leq 20$ (attained at $n=60, 84$) and $R(33; 8) < 7.6$.
If $20 \leq n \leq 32$, then $u \leq 11$ and $R(20;8) <7.8$. The
values of $n= 16, 17, 19$ are excluded by previous considerations:
the remaining values  $n=15$ or $18$ have $u=5$ so that $R(15; 8)
<6.2$. Thus Case I  has been completed simply by additive sieving
with $\rho=1/3$.

 \bigskip
\noindent {\bf Case II:}  $q= 5$ or $7$.

 \noindent This follows broadly the same pattern as Case I, except
 that, because $2^{8/3} > 6.34$,  the expression
 $R_2(n)$ is   useless  when $q=5$ and ineffective when $q=7$.
 We therefore proceed as follows. Suppose $\n > q^2$ so that $s \geq 3$.
  Suppose first that also  $\omega(Q) \geq 49$. By Lemma $\ref{1/6bound}$ and the fact that
 $\frac{n}{2}-\frac{n}{6}=\frac{n}{3}$, we obtain as an   alternative to
 $(\ref{weerhoeq})$
 \begin{equation}\label{weerhoeq1}
R_3(n)=R_3(n;q,s)\ =\ \left\{2^{(2/3)
n+1}\frac{1}{(q-1)^{1/6}}\left(\frac{\frac{4n}{3s}
-1}{{1-\frac{4n}{3sq^{s}}}} +2\right)\right\}^{3/n}.
\end{equation}
 Here the starter form, derived from $(\ref{weerhoeq1})$ using $ s \geq 3$
and $n < q^s$, is
\begin{equation}\label{preweerho1}
\bar{R}_3(n)= \bar{R}_3(n;q)\ =\ \left\{c2^{2\rho
n+1}\frac{(4n+1)}{5(q-1)^{1/6}}\right\}^{3/n}.
\end{equation}
Now $\bar{R}_3(58;5)< 4.998$ and $\bar{R}_3(16;7) < 6.98$.
Summarising, whenever  $\w(Q) \geq 49$, we have shown that
necessarily $n \leq 57$ ($q=5$) and $ n \leq 15$ ($q=7$).  But,
easily, if $n \leq 57$ (say), then $\w(Q)<49$.

Hence we may suppose  that $\w(Q)\leq 48$. Since $s \geq 3$ the
appropriate starter form for $R(n)$ itself (in place of
$(\ref{starterR})$) is
\[ \bar{R}(n)  =\bar{R}(n;5,u): = \ \{2^{(2/3)
n+u+1}(4n+1)/5\}^{2/n}.\]

 For the rest,  we focus almost exclusively on  the
more delicate case when $q=5$.
 Then with   $\bar{R}(113; 5,48) <
4.98$. So assume $n \leq 112$  in which case since $Q\leq
(5^n-1)/4$, necessarily $\w(Q) \leq 44$. Moreover, since $R(104;5)
<4.99$ (with  $\rho =1/3$ and $s=2$), we can suppose that $n \leq
103$. Indeed, by repetition of this argument using $R(n;5,u)$ and
smaller values of $u$, we conclude that we can suppose $n \leq 84$.

The next stage (with $n\leq 84$) is to calculate the true value of
$\w(Q)$ and use  $R(n)$ (still with $\rho=1/3$ and $s=2$).  We find
that $R(44;5)<4.98$ so that we can assume $n\leq 43$.  But then
$\w(Q) \leq  11$ and, with $u=11$, $R(33;5)<4.97$. Next, $n \leq 29$
and $\w(Q)\leq 8$ and with $u=8$, $R(26; 5)<4.98$.  Further, with
$u=4$, $R(21;5)<4.14$. The values of $n$ that remain belong to the
the set $\{9, 12,18,24\}$. When the same exercise is applied to the
field with $q=7$, the only outstanding degree is $n=9$.  We tabulate
the outcome  of applying Proposition $\ref{rhoform}$ in full, in one
 case using the  form $(\ref{Rform})$ for $R(n)$.

\bigskip
\bigskip
\begin{tabular}{ccccccccc}
  \hline
   $q$ & $n$&$s$&$\rho$&$Q$&$u$&$t$&$\delta$&$R(n)$ \\ \hline \hline
 $5$&$9$&$6$&$2/9$&$19\cdot 31 \cdot 829$&$3$&$0$&$1$&$4.49$\\ \hline
       $5$&$18$&$6$&$2/9$&$3^3\cdot7\cdot19\cdot 31\cdot 829\cdot 5167$&$4$&$0$&$1$&$ 3.85$\\
\hline
$5$&$24$&$2$&$1/6$&$2\cdot3^2\cdot7\cdot13\cdot31\cdot313\cdot601\cdot390001$&$8$&$0$&$1$&$3.91$\\
\hline

   $7$&$9$&$3$&$1/3$&$3\cdot 19\cdot 37 \cdot 1063$&$1$&$3$&$0.919$&$4.82$\\ \hline

\end{tabular}

\bigskip
\bigskip
 For the pair $(5,12)$,  Proposition $\ref{rhoform}$ fails:  in that
 case we found the explicit PFF polynomial
 $x^{12}+x^{11}+x^3-x^2-2x-2$.

\bigskip
As a consequence, Proposition $\ref{weerho}$ is established.

\end{proof}

\bigskip
\bigskip

\section{Very small fields}\label{vsmall}

In this section, we study the smallest fields $\mathbb{F}_q$ when $2
\leq q \leq 4$.  For these it is imperative to use a smaller value
of $\rho$ than provided by Lemma $\ref{rho}$.  Variations of Lemma
$\ref{1/6bound}$ are also invoked where appropriate. Further,
 more attention has to be paid than heretofore when $\n<n$: in
particular the refined forms of Lemma $\ref{Rform}$ will be called
on  to resolve some smaller values.

\begin{lemma}[\cite{CoHu03}]\label{rhowee}
Assume that $n>4$ ($p \nmid n$).   Then the following hold.
\begin{itemize}
\item[(i)] Suppose $q=4$.  Then $\rho(4,9) =1/3; \ \rho(4,45)
=11/45;\ $ otherwise $\rho(4,n) \leq 1/5$.

\item[(ii)] Suppose $q=3$.  Then $\rho(3,16) =5/16;\ $ otherwise
$\rho(3,n) \leq 1/4$.

\item[(iii)] Suppose $q=2$.  Then $\rho(2,5) =1/5; \ \rho(2,9)
=2/9;\ \rho(2,21) =4/21;\ $otherwise $\rho(2,n) \leq 1/6$.
\end{itemize}
\end{lemma}

\subsection{The field $\mathbb{F}_4$} \label{four}

Here $n^*= n$ if and only if $n$ is odd, whereas $Q$, a divisor of
$(4^n-1)$, is always odd.

\medskip
\begin{prop}\label{4result}
Suppose $q = 4$ and $n \neq 3$. Then $(q,n)$ is a PFF pair.
\end{prop}

\begin{proof}
For the main working suppose $n^* >4$ and $s>1$. For $n$ odd, by
Lemma $\ref{rhowee}$,  $\rho(n) \leq 1/5$, except when $n=9$ ($\rho=
1/3$) or $n=45\ (\rho=11/45)$. When $n$ is even, $\rho \leq 1/6$
(with equality when $n=18$). Further, $s=2$ when $n^*$ divides $15$;
$s=3$ when $n^*$ divides $63$; otherwise $s \geq 4$.

\bigskip
Start from the  sufficient condition $(\ref{ineq1})$ with $R(n)$
given by $(\ref{Rformadd})$ and $u=\w(Q)$.

First suppose $\n=15$ (the only situation in which Lemma
$\ref{rhoincr}$ does \emph{not} apply); thus $\rho=1/5$. Since  the
expression
\begin{equation}\label{rfrac}
E=\frac{\frac{2(\n-\rho n)}{s} -1}{{1-\frac{2(\n-\rho n)}{sq^s}}}+2
\end{equation}
 in the \emph{refined form} of $(\ref{Rformadd})$ here is equal
to $46$ and the (crude) bound
 $W(Q)<2.9 \cdot 4^{n/4}$ holds (by $\ref{Wbound}$), it follows
 that inequality $(\ref{ineq1})$ certainly holds whenever
 \[ 4> (2.9\cdot2\cdot 46)^{20/n}=266.8^{20/n},\]
 and this is satisfied when $n \geq 120$.  Thus, when $\n=15$ it can
 be assumed that $n \leq 60$.
 Now suppose that $\n \neq 15$.  In order to
construct a suitable starter function for larger values of $n$, by
Lemma $\ref{rhoincr}$ replace $\rho$ by a larger value (such as
$1/5$ or $1/6$).
 To $W(Q)=2^\omega$, again apply the  bound
$W(h) < 2.9h^{1/4}$ (Lemma $\ref{Wbound})$. Using $n<q^s$ and $s
-2(1-\rho)<s$, we see that $4>R_3(n)$ suffices, where
\[
R_3(n)\ =R_3(n;s,\rho)\
\left\{5.8\left(\frac{2(1-\rho)n}{s-2(1-\rho)}+1\right)\right\}^{\frac{4}{4(1-4\rho)}}
\]
with the appropriate larger value of $\rho$. Here $R_3(n)$ decreases
as a function of $s$ and decreases as a function of $n$.

If   $n \ (\neq 45) $ is odd and $s\geq 4$, then $\rho \leq 1/5$
 and $R_3(n;4,1/5)<83.4$.  If $n$ is even then $\rho \leq 1/6$ and $R_3(n;2,1/6)<65.5$.
Since $ s\geq 4$ whenever $\n >63$, it follows that for a putative
exception $n$ to  Proposition $\ref{4result}$ we may assume $n \leq
83$; indeed, $n \leq 62$ for $n$ even.

For these remaining possibilities (including those with $\n=15)$, we
evaluate $R(n)$ given by (the refined form of) $(\ref{Rformadd})$
with precise values for $s, \rho$ and $u=\w(Q)$: if it is  less than
$q=4$ there \emph{does} exist a PFF polynomial for that value of
$n$. To this end factorise $x^{\n}-1$ over $\mathbb{F}_4$ and $Q$.
For larger values of $n$ and those for which $\n$ is prime (in which
case $\rho=1/n$), comfortably $R(n) <4$.  We tabulate the outcome in
the more delicate cases with $n \geq 10$: in particular, the column
headed $R$ lists $R(n)$ truncated to three decimal places.

\bigskip\begin{tabular}{c c ccc || ccccc}
\hline
  $n $ & $s$ &$\rho$& $u$ & $R$ & $n $ & $s$ &$\rho$& $u$ & $R$  \\ \hline \hline
  $45$&$6$&$11/45$&$11$&$3.187$&$21$&$3$&$1/7$&$6$&$3.063$\\ \hline
 $36$&$3$&$1/12$&$12$&$2.277$&$20$&$2$&$1/20$&$7$&$2.392$\\ \hline
  $35$&$6$&$1/7$&$9$&$2.532$&$18$&$3$&$1/6$&$8$&$3.815$\\ \hline
   $33$&$5$&$1/11$&$8$&$2.195$&$15$&$2$&$1/5$&$6$&${\bf 5.539}$\\ \hline
    $30$&$2$&$1/10$&$11$&$2.965$&$14$&$3$&$1/14$&$6$&$3.085$\\ \hline
     $27$&$9$&$5/27$&$6$&$2.729$&$11$&$5$&$1/11$&$4$&$3.238$\\ \hline
      $25$&$10$&$1/5$&$4$&$3.238$&$10$&$2$&$1/10$&$5$&${\bf 4.337}$\\ \hline

\end{tabular}

\bigskip

We conclude that if there is no PFF polynomial of degree $n$, then
$n \in \{15,10,9,7,5\}$.  For the values $n=10, 7$, using also
multiplicative sieving yields  the result. Specifically, suppose
$n=10$.  Then $Q=3^2\cdot5^2 \cdot11\cdot 31\cdot 41$, which has 5
prime factors. In $(\ref{Rform})$, take $u=1, \ t=4$.  Then
$\delta>0.6524$ which yields $R(10)<3.73 <4$. For $n=7$,
$Q=3^2\cdot43 \cdot127$.  In this case, take $u=1, \ t=2$, so that $
\delta > 0.9688$ and $R(7)< 3.93<4$.

Finally,  we exhibit explicit  PFF polynomials  for the remaining
degrees (including $n=6$, held over from Proposition $\ref{n=q-1}$).
For these, we use $\mathbb{F}_4=\mathbb{F}_2(u)$, where $u^2+u+1=0$.

\bigskip
$
\begin{tabular}{c|c}
  \hline
   $n$ & PFF polynomial \\ \hline \hline
     $15$ & {\tiny$ x^{15}+x^{14}+(u+1)x^{12}+(u+1)x^{10}+x^9+x^8+x^7+ux^6+ux^5+ux^4+x^2+ux+u+1$}\\ \hline

    $9$ &  $x^9+(u+1)x^8+ux^7+(u+1)x^6+ux^5+ux^3+(u+1)x+u$\\ \hline
$6$ & $x^6+ux^5+(u+1)x^4+(u+1)x^3+x+u+1$ \\ \hline
  $5$ & $x^5+ux^4+ux^3+x+u+1 $\\ \hline
\end{tabular}
$

\bigskip
\bigskip

\end{proof}

\subsection{The ternary field $\mathbb{F}_3$}\label{three}

For the main part, again suppose $n^*>4$ and $s\geq 2$. Here  any
version of Lemma $\ref{Wbound}$ valid for all integers is
inadequate:  the following numerical bound for  large integers will
be needed.
\begin{lemma} \label{hboundq=3}
Suppose $h$ is indivisible by $3$ and $\w(h) \geq 52$. Then
\[ W(h) < h^{4/25}.\]
\end{lemma}

\medskip
\begin{prop}\label{3result}
Suppose $q = 3$ and $n \neq 3$. Then $(q,n)$ is a PFF pair.
\end{prop}

\begin{proof}
By Lemma $\ref{rhowee}$, if $n^*=n$ (equivalent to $3 \nmid n$),
then   $\rho(n) \leq 1/4$, except when $n=16$.  If, on the other
hand, $3|n$, then evidently, $\rho(n) \leq 5/48$; indeed $\rho(n)
\leq 1/12 $ whenever $n > 48$.

\bigskip
Again, start from the  sufficient condition $(\ref{ineq1})$ with
$R(n)$ given by $(\ref{Rformadd})$ and $u=\w(Q)$.

Suppose $3|n$ (i.e.,  $n>\n$) with $n>48$ so that $\rho \leq 1/12$:
in this situation Lemma $\ref{Wbound}$ suffices. Since $\n<n, \n\leq
q^s-1, s\geq 2$ and $\rho
>0$, then $E$ (given by $(\ref{rfrac})$) satisfies

\[E<\frac{\frac{2n}{s} -1}{{1-\frac{2\n}{sq^s}}}+2\leq 9n-7.\]
By Lemma $\ref{Wbound}$ and the fact that $\rho \leq 1/12$, it
suffices that $3^{1/2}/2^{1/3} > (9.8\cdot(9n-7))^{2/n}$, which
holds whenever $n\geq 54$. Hence we may assume $n \leq 51$ when
$3|n$.

Now suppose $3 \nmid n$ (so that $\n=n$). With Lemma
$\ref{hboundq=3}$ in view, suppose $\w(Q) \geq 52$ so that certainly
$\rho \leq 1/4$ and $s \geq 4$.  Since $n=\n \neq 8$, in $R(n)$
replace $\rho$ by $1/4$, as we may by Lemma $\ref{rhoincr}$.  From
$(\ref{Rformadd})$ and Lemma $\ref{hboundq=3}$, we derive the
sufficient condition
\[3^{17/25}/2> \left(\frac{2(3n+2)}{5}\right)^{2/n},\]
which holds whenever $n \geq 205$ and therefore whenever $\w(Q) \geq
52$.

Continue to suppose $3 \nmid n$  with $n \geq 55$ and $n\neq 80$ (so
that $\rho \leq 1/4$ and $s \geq 5$) but assume now that $W(Q) \leq
51$. We introduce a multiplicative aspect to the sieve by invoking
$R(n)$ as in $(\ref{Rform})$.
  To show that that $R(n)$ is increasing with $\rho$ analogously
to Lemma $\ref{rhoincr}$ consider
\begin{equation}\label{diff1}\tau(\rho)= \log\left[2^{2\rho
n}\left(\frac{2(1-\rho)(n/s)+t-1}{\delta-(2(1-\rho)n/sq^s)}\right)\right],
\end{equation}
with $q=3$. Here we suppose $\delta$ is bounded below by $0.42$, an
assumption that will be realised in applications. (In the first
place, since $\rho >0$ and $s \geq 5$, this guarantees that
$\delta-(2(1-\rho)/s)$ and so $\delta-(2(1-\rho)n/sq^s)$ are
positive.)
 For fixed
$s$, differentiate $\tau(\rho)$ to obtain
\begin{equation}\label{diff2}\tau^\prime(\rho)=2n\log 2-
\frac{1}{(1-\rho)+(s(t-1)/2n)}-\frac{1}{(\delta sq^s/2n)-(1-\rho)},
\end{equation} with $q=3$. Since $0< \rho \leq 1/4, n<3^s$, $s\geq 5$
and $\delta \geq 0.42$ it follows that $\tau^\prime(\rho)\geq 2n\log
2-4/3-20 =2n\log 2-64/3 $ which is positive because $n \geq 16$.

Granted that $\delta \geq 0.42$ it can be concluded that, for a
given $n$ and $t$, $\tau(n)$ and so $R(n)$ are maximised when $s=5$
and $\rho=1/4$.
 This yields the condition $3 > R_4(n)$, where
 \begin{equation}
 \label{eq4q=3}
 R_4(n)\ =\ 2\left( 2^{1+u}\left(\frac{3n+10(t-1)}{10\delta-3}+2\right)\right)^{2/n}.
\end{equation}
with $t$ denoting the number of sieving primes and $u$ those of the
multiplicative core $m_0 \ (| Q)$.
 To use  $(\ref{eq4q=3})$, let the least $u=6$  primes in $Q$ contribute to the core $m_0$.
 Then $t\leq 45$ is the number of sieving primes and $\delta \geq
\frac{1}{19}+\frac{1}{23}++\cdots+\frac{1}{239}= 0.42734\ldots\ $.
   Since  $R_4(55.4)<3$ there  exists a PFF polynomial of degree $n$ whenever  $n \geq 55$ ($n \neq 80$).

Summarising, whether or not $3|n$, it remains to consider values of
$n \leq 53$ and $n=80$. One could apply further general applications
of the sieve to some  effect but instead we simply calculate $R(n)$
given by (the refined form of) $(\ref{Rformadd})$.  In the table,
the column headed $R$ gives its value truncated to three decimal
places. Only those degrees which produced a value of $R(n)$
exceeding $2.2$ are listed: none of these has $\n<n$.

 \bigskip\begin{tabular}{ccccc||ccccc}
\hline
  $n $&$s$&$\rho $&$u$ & $R$&  $n $&$s$&$\rho$&$u$ & $R$\\ \hline \hline
$52$&$6$&$11/52$ & $6$ & $2.390$ & $14$ & $6$ & $1/7$ &$3$&$2.780$
\\\hline
$44$&$10$&$7/44$ & $8$ & $2.245$ & $13$ & $3$ & $1/13$ &$1$&$2.243$
\\\hline
$32$&$8$&$7/32$ & $6$ & $2.811$ & $11$ & $5$ & $1/11$
&$2$&$2.520$\\\hline
 $28$&$6$&$3/28$ & $6$ & $2.234$ & $10$ & $4$ & $1/5$ &$3$&${\bf
4.208}$\\\hline

$22$&$5$&$1/11$ & $5$ & $2.298$ &$8$ & $2$ & $1/4$ &$3$&${\bf
8.122}$\\ \hline
 $20$&$4$&$3/20$ & $5$ & $2.903$ &  $7$ & $6$ & $1/7$ &$1$&${\bf 3.023}$\\ \hline

 $16$ & $4$ & $5/16$ &$4$&${\bf 5.085}$&$5$ & $4$ & $1/5$
&$1$&${\bf 4.720}$\\ \hline
\end{tabular}

\bigskip
\bigskip
To supplement this table note that when  $n=7$ we can successfully
use $(\ref{Rform})$ by sieving also with the single prime divisor of
$Q=1093$: this yields $R(7)<2.694<3$. Including  cases held over
from Proposition $\ref{small n^*}$, this leaves $ n \in
\{16,12,10,8,6,5,3\}$ for which we obtain a PFF polynomial in every
case by direct verification of the properties. In fact when $n=3$
there is only one pair of PFF polynomials.

\bigskip
\bigskip\begin{tabular}{c|c}
  \hline
   $n$ & PFF polynomial \\ \hline \hline
   $16$ & $x^{16}-x^{15}-x^6+x-x-1$\\ \hline
  $12$ & $x^{12}+x^{11}+x^3+x^2+x-1$\\ \hline
  $10$ & $x^{10}+x^9+x^7+x^3-x-1$ \\ \hline
  $8$ & $x^8+x^7+x^4-x^3-x^2+x-1$ \\ \hline
  $6$& $x^6+x^5+x^3+x^2+x-1$ \\ \hline
  $5$&$x^5+x^4-x+1$\\ \hline
        $3$& $x^3+x^2-x+1$\\ \hline
\end{tabular}

\end{proof}

We remark that we incorporated multiplicative sieving as a device to
treat general values of $n\geq 55$ (with $3\nmid n$) in Proposition
$\ref{3result}$. Nevertheless, it is likely that for any
\emph{specific} value of $n \geq 55$ additive sieving using
$(\ref{Rformadd})$ would be sufficient. A similar remark would apply
to the proof of Proposition $\ref{2result}$ below.

\bigskip
\bigskip

\subsection{The binary field $\mathbb{F}_2$.}\label{two}
   A suitable numerical result on $W(h)$ here
 is the following.
 \begin{lemma}\label{2W}
Suppose the \emph{odd} integer $h$ is such that $\omega(h) \geq
175$.  Then $W(h) < h^{3/25}$.
\end{lemma}

\begin{prop}\label{2result}
Suppose $q = 2$ and $n \neq 3,4$. Then $(q,n)$ is a PFF pair.
\end{prop}
\begin{proof}
The cases $(2,n),\ n=6,12,24$ have been held over from  Proposition
$\ref{small n^*}$.  Otherwise, suppose that $n^*>4$, so that $s \geq
3$.  Here $Q=2^n-1$. By Lemma $\ref{rhowee}$, if $n$ is odd (i.e.,
$\n=n$) and $n>21$, then $\rho\leq 1/6$.  If $n$ is even and then
$\rho \leq 1/12$ unless $\n=5,9$ or $21$; indeed, if $n
>42$, then $\rho \leq 1/12$.

Suppose first that  $n >42$ is \emph{even}, so that $\rho \leq
1/12$. In Proposition $\ref{rhoform}$, since $0<\rho \leq 1/12$ then
$R(n)$ given by the refined form of  $(\ref{Rformadd})$ satisfies
\begin{equation}\label{neven}
R(n) <
\left\{W(Q)\cdot2^{\frac{n}{6}+1}\cdot\left(\frac{\frac{2\n}{s}-1}{1-\frac{2\n}{sq^s}}+2\right)\right\}^{2/n}
\leq
\left\{W(Q)\cdot2^{\frac{n}{6}+1}\cdot\frac{(4n-2)}{5}\right\}^{2/n}.
\end{equation}
Since $W(Q)< 2.9\cdot 2^{n/4}$, it follows that there exists a PFF
polynomial of degree $n$ whenever $2>\left\{5.8
\cdot\frac{(4n-2)}{5}\right\}^{12/n}$ and so certainly when $n \geq
108$.

 The general argument with $n$ even is taken somewhat further.
  Suppose $n \leq 106$.  By calculation, $\w(Q)\leq 21$.
  Substituting $W(Q)=21$ in $(\ref{neven})$ we find that $R(86) <2$;
  hence we may suppose $n \leq 84$.  Indeed, by  repetition of this strategy we conclude there exists a
  PFF polynomial of degree $n$ whenever $n>64$.

Now suppose $n\ (>64)$ is \emph{odd} so that $\n=n$, $\rho \leq 1/6$
and $s \geq 7$. By Lemma $\ref{rhoincr}$ we can replace $\rho$ by
$1/6$ in $R(n)$ given by $(\ref{Rformadd})$. In order to apply Lemma
$\ref{2W}$ suppose (temporarily) that additionally $\w(Q) \geq 175$.
Since $1/2-1/3-3/25=7/150$, $n<q^s$ and $s\geq 7$ we deduce that
there is a PFF polynomial of degree $n$ whenever
 \[2 \ >\  R(n)=\left(\frac{2(5n+11)}{16}\right)^{150/7n} \]
and so whenever $n \geq 139$.  Easily, this  is implied by $\w(Q)
\geq 175$.

Accordingly, we can now suppose $\w(Q) \leq 174$.  Introduce a
multiplicative dimension to the sieve by applying the criterion of
Lemma $\ref{rhoform}$ with $R(n)$ given by $(\ref{Rform})$. By
$(\ref{diff1})$ with $q=2$ and provided $\delta >0.42$, $\tau(n)$ is
increasing for $0 < \rho \leq 1/6$, since $\tau^\prime(\rho)\geq
2n\log 2-7/6-100/47 \geq 2n\log 2-4 $ is positive.
 Hence in $R(n)$ we may replace $\rho$ by $1/6$ and $s$
by $7$, to obtain the sufficient condition
\[ 2\ >\  R_5(n):= \left\{2^{u+1}\left( \frac{5n+21(t-1)}{21\delta-5}+2
\right)\right\}^{6/n},\] provided $\delta >0.42$, where $u$ is the
number of prime integers in the multiplicative core.

First take $u=13$ so that $t \leq 161$. Then $\delta > 0.4354$ and
$R_5(144) < 2$. Hence we can suppose $n \leq 143$.  This implies
$\omega \leq 27$. Thus $u+t \leq 28$.  Repeat the above process with
$u=4, \ t \leq 23$ and so $\delta >  0.4353$.  Then $R(77)< 2$ and
we can suppose $n \leq 75$.   Then $\w(Q) \leq 16$. Repeat once more
with $u=3$, $\delta > 0.4787$ to yield $R(66)< 2$.

Consequently, for the last stage, whether $n$ is even or odd, assume
$n \leq 65$.  As for Lemma $\ref{3result}$, simply calculate $R(n)$
given by (the refined form of) $(\ref{Rformadd})$. The table lists
the outcome for values of $n$ with $13 \leq n \leq 65$ which
produced a value of $R(n)$ exceeding $1.8$. Also included is $n=24$
with $s=2$, a case held over from  Proposition $\ref{small n^*}$.

 \bigskip\begin{tabular}{ccccc||ccccc}
\hline
  $n $&$s$&$\rho $&$u$ & $R$&  $n $&$s$&$\rho$&$u$ & $R$\\ \hline \hline
$45$&$12$&$2/15$ & $6$ & $1.963$ & $24$ & $2$ & $1/24$ &$6$&$1.887$
\\\hline
$42$&$6$&$2/21$ & $6$ & $1.801$ & $22$ & $10$ & $1/22$ &$4$&$1.717$
\\\hline
$36$&$6$&$1/12$ & $8$ & $1.895$ & $21$ & $6$ & $4/21$ &$3$&${\bf
2.662}$\\\hline
 $35$&$12$&$4/35$ & $4$ & $1.856$ & $20$ & $4$ & $1/21$ &$5$&$1.941$\\\hline

$30$&$4$&$1/15$ & $6$ & $1.953$ &$18$ & $6$ & $1/9$ &$4$&${\bf 2.290}$\\
\hline
 $28$&$3$&$1/28$ & $6$ & $1.811$ &  $15$ & $4$ & $2/15$ &$3$&${\bf 2.892}$\\ \hline

 $27$ & $3$ & $1/9$ &$3$&$1.839$&$14$ & $3$ & $1/14$
&$3$&${\bf 2.438}$\\ \hline

 $25$&$20$&$2/25$ & $3$ & $1.714$ & $13$ &
$12$ & $1/13$ &$1$&$1.814$\\\hline
\end{tabular}

\bigskip

Beyond this table, degrees   $n=11,18$ and $21$ can be treated
theoretically. For $n=11$ use $(\ref{Rform})$ by sieving also with
the two prime divisors of $Q=23\cdot89$. This yields $R(11)=
1.968\ldots <2$. Similarly, when $n=18$, sieve also with the $4$
prime divisors of $Q=3^3\cdot 7\cdot19\cdot 73$. This yields $R(18)=
1.980\ldots <2$. Finally when  $n=21$, for this occasion only,
modify the key strategy for the additive sieve as follows. Over
$\mathbb{F}_2$, $x^{21}-1=P_1\cdot P_2\cdot P_{31}\cdot P_{32}\cdot
P_{61}\cdot P_{62}$, where the $P$'s are distinct irreducible
polynomials of degree indicated by the first subscript. For the
sieve take the ``core" to be $P_1P_2$ and the sieving irreducibles
to be those of degrees $3$ and $6$ together with the three  prime
factors of $Q=7^2\cdot127\cdot337$. The crucial denominator $\delta
-4/2^3-4/2^6=0.2838\ldots $ and $R(21)=1.963\ldots<2$

For degrees $\n \in \{15,14,12,10,9,8,7,6,5\}$ (including some held
over from Proposition $\ref{small n^*}$),  we obtain a PFF example
in every case. We remark that, for $n=6,5$ there is, in each case, a
single pair of reciprocal PFF polynomials.

\bigskip
\bigskip\begin{tabular}{c|c}
\hline
   $n$ & PFF polynomial\\ \hline \hline
  $15$ & $x^{15}+x^{14}+x^4+x+1$\\ \hline
  $14$ & $x^{14}+x^{13}+x^9+x^4+x^2+x+1$ \\ \hline
  $12$ & $x^{12}+x^{11}+x^9+x^4+x^3+x+1$ \\ \hline
  $10$& $x^{10}+x^9+x^4+x+1$\\ \hline
  $9$&$x^9+x^8+x^5+x^4+x^3+x+1$\\ \hline
        $8$& $x^8+x^2-x+1$\\ \hline
  $7$&$x^7+x^4+x^3+x+1$
  \\ \hline
        $6$& $x^6+x^5+x^2+x+1$\\ \hline
  $5$&$x^5+x^4+x^2+x+1$\\ \hline
\end{tabular}

\bigskip
\bigskip

\end{proof}

 {\bf Acknowledgement} The second author is supported by a Royal Society Dorothy
Hodgkin Fellowship.

 \vspace{5mm}
\noindent{Stephen D. Cohen \\ Department of Mathematics \\
University of Glasgow \\ Glasgow, G12 8QW, UK \\ Email:
sdc@maths.gla.ac.uk}

\vspace{5mm}

\noindent{Sophie Huczynska \\ School of Mathematics and Statistics \\
University of St Andrews \\ St Andrews, Fife, KY15 7NA, UK \\
Email: sophieh@mcs.st-and.ac.uk}

\end{document}